\documentclass[11pt,letterpaper]{article}

\usepackage{color}

\usepackage{epsfig,lscape,rotate}
\usepackage{amsmath}
\usepackage{lscape}
\usepackage{amssymb}
\usepackage{amsfonts}
\usepackage{euscript}
\usepackage{mathrsfs}
\usepackage{natbib}
\usepackage{multirow}

\usepackage{fullpage}
\newtheorem{thm}{Theorem}

\newtheorem{lema}{Lemma}

\newlength{\defaultbaselineskip}
\newcommand{\Lvb}{L\'{e}vy }

\newcommand{\Ito}{It\^o }
\newcommand{\cadlag}{c\`{a}dl\`{a}g }
\newcommand{\Holder}{H\"{o}lder}

\newcommand{\refeq}[1]{(\ref{#1})}
\newcommand{\qedd}{$\hfill\Box$}

\newcommand{\MR}{$\color{red}\spadesuit$}
\DeclareMathOperator{\R}{{\mathbb R}}

\DeclareMathOperator{\N}{{\mathbb N}}

\providecommand{\eps}{\varepsilon}
\renewcommand{\phi}{\varphi}
\renewcommand{\theta}{\vartheta}

\renewcommand{\cdot}{{\scriptstyle \bullet} }
\providecommand{\abs}[1]{\lvert #1 \rvert}
\providecommand{\norm}[1]{\lVert #1 \rVert}

\renewcommand{\le}{\leqslant}
\renewcommand{\ge}{\geqslant}

\begin{document}

\begin{titlepage}
\title{Nonparametric test for a constant beta between \Ito semi-martingales based on high-frequency data\footnote{Todorov's work was partially supported by NSF grant SES-0957330, Rei{\ss}'s work was partially supported by the DFG via SFB 649 {\it Economic Risk}.}}
\date{\today}
\vspace{1.0cm}
\author{\\\mbox{}\\\mbox{} Markus Rei{\ss}\thanks{Institut f\"{u}r Mathematik, Humboldt-Universit\"{a}t zu Berlin, Unter den Linden 6, 10099 Berlin, Germany;
email: mreiss@math.hu-berlin.de.}~~and~~Viktor Todorov\thanks{\noindent Department
of  Finance, Northwestern University, Evanston, IL 60208-2001, email:
v-todorov@northwestern.edu.}~~and~~George Tauchen\thanks{Department of Economics, Duke University, Durham, NC 27708; e-mail: george.tauchen@duke.edu.}\\\mbox{}}
\end{titlepage}

\maketitle

\begin{abstract}
\noindent We derive a nonparametric test for constant  beta over a fixed time interval from high-frequency observations of a bivariate \Ito semimartingale. Beta is defined as the ratio of the spot continuous covariation between an asset and a risk factor and the spot continuous variation of the latter. The test is based on the asymptotic behavior of the covariation between the risk factor and an estimate of the residual component of the asset, that is orthogonal (in martingale sense) to the risk factor, over blocks with asymptotically shrinking time span. Rate optimality of the test over smoothness classes is derived.
\end{abstract}

\vspace{1.5cm}



\renewcommand{\arraystretch}{1.0}\setlength{\baselineskip}{6.5mm}
\providecommand{\MR}{$\clubsuit$}

\bigskip \newpage

\section{Introduction}

In this paper we develop a test for the time-variation of the process $\frac{d\langle X^c, Y^c\rangle_t}{ d\langle X^c, X^c\rangle_t }$ over a fixed interval of time. Here $X$ and $Y$ are two semimartingales, $X^c$ and $Y^c$ denote their continuous components, and the angle bracket denotes the predictable component of the quadratic (co)variation, see e.g., \cite{JS}. Our analysis applies to bivariate \Ito semimartingales defined on a filtered probability space $\left(\Omega,\mathcal{F},(\mathcal{F}_t)_{t\geq 0}, \mathbb{P}\right)$ and having representation of the following form
\begin{equation}\label{eq:XY}
X_t = X_0+F_t^X+\int_0^t\sigma_sdW_s,~~Y_t = Y_0+F_t^Y+\int_0^t\beta_s\sigma_sdW_s+\int_0^t\widetilde{\sigma}_sd\widetilde{W}_s,
\end{equation}
where $X_0$ and $Y_0$ are $\mathcal{F}_0$-measurable random variables, $F^X$ and $F^Y$ are finite variation processes (containing both continuous and jump parts),  $W$ and $\widetilde{W}$ are two independent Brownian motions, $\beta$, $\sigma$ and $\widetilde{\sigma}$ are stochastic processes with \cadlag paths, exact assumptions being provided in the next section. In the setting of \refeq{eq:XY} the continuous quadratic covariation is absolutely continuous with respect to time and $\beta_t \equiv  \frac{d\langle X^c, Y^c\rangle_t}{ d\langle X^c, X^c\rangle_t }$, and hence our interest in this paper is in testing whether the process $\beta$ remains constant over a given time interval. The key motivating example for this problem comes from finance where $X$ plays the role of a risk factor and $Y$ of an asset. $\beta$ in this case measures the exposure of the asset to the risk factor, and constancy of $\beta$ plays a central role in testing the validity of the asset pricing model. Indeed, the time-variation in beta can generate an excess return in the asset, above what is implied by the model, and hence lead to its rejection. Testing whether $\beta$ is constant on a given interval helps further decide on the time window for recovering the beta process from the data.

The asymptotic analysis in the paper is based on discrete equidistant observations of the bivariate \Ito semimartingale on a fixed interval of time with mesh of the observation grid shrinking to zero. Our focus on the sensitivity of the continuous martingale part of $Y$ towards that of $X$ is similar to \cite{BNS2004Ecma}, \cite{ABDG2006}, \cite{TB} and \cite{GM}.

We construct our test as follows. We first form a ``pooled'' estimate of beta as the ratio of estimates over the fixed interval $[0, T]$ of the continuous covariation $\langle X^c, Y^c\rangle_T$ and the continuous variation $\langle X^c, X^c\rangle_T$.  This estimator is consistent for the constant beta and asymptotically mixed normal under the null and it converges to a volatility weighted average of the time-varying beta under the alternative. Using this ``pooled'' beta estimator, we then separate, under the null hypothesis of constant beta, the residual component of the process $Y$ which is orthogonal in the continuous martingale sense to the process $X$. That is, we estimate, under the null, the part of $Y$ that has zero continuous quadratic covariation with $X$. Since the ``pooled'' beta estimates the true beta process only under the null, the above estimate of the residual component is asymptotically orthogonal to $X$ only when beta is constant.

Our test is formed by splitting the data into blocks of decreasing length and forming test statistics for constant beta on each of the blocks. This is similar to block-based estimation of volatility functionals in high-frequency setting developed in \cite{JLMPV}, \cite{MZ08} and \cite{JR}. The test statistics on the blocks are based on the different asymptotic behavior under the null and alternative of our estimate of the residual component defined in the previous paragraph. Our test is then formed by summing the test statistics over the blocks and appropriately scaling the resulting sum. The test is asymptotically standard normal under the null and after scaling it down it converges to a volatility weighted measure of dispersion of the beta around its volatility weighted average on the fixed time interval.

The asymptotic behavior of our statistic has several distinctive features compared with block-based volatility functional estimators considered in \cite{MZ08} and \cite{JR}. To achieve non-degenerate limits under the null of constant beta, unlike \cite{JR}, we need to scale up appropriately the local block variance-covariance estimates. As a result, unlike \cite{JR}, the limiting distribution of our statistic is not determined from the first-order expansion of the nonlinear function of the block volatility estimates around the function evaluated at the true (and observed) stochastic variance-covariance matrix. We further extend the analysis in \cite{JR} by considering functions of volatility which are not bounded around zero. Finally, unlike earlier work, our statistic is constructed as a nonlinear function of adjacent volatility block estimators. This makes the effect of biases arising from the local volatility estimation negligible and in particular it circumvents the need to do any bias correction which from a practical point of view is very desirable.

Turning to the testing problem, our test has three distinctive features. First, the test is pathwise in the sense that it tests whether beta is constant or not on the observed path. Hence the analysis here is based on in-fill asymptotics and it requires neither assumptions regarding the sources  of the variation in beta nor stationarity and ergodicity conditions. Second, our test statistic is of self-normalizing type (see \cite{SN}) and hence its limiting distribution under the null is pivotal, i.e., it does not depend on ``nuisance parameters'' like the stochastic volatilities of the two processes. Finally, we can show that our test is asymptotically optimal for local nonparametric alternatives $\beta_t$ that are $\alpha$-H\"older regular. The separation rate of a weighted $L^2$-distance between hypothesis and alternative is $n^{-2\alpha/(4\alpha+1)}$, for which a minimax lower bound proves its optimality. This analysis also provides a rationale for selecting the block size, depending on which kind of alternatives we would like to discriminate. Let us also remark that a simple test based on the difference of a nonparametric estimator of $\beta_t$ from a constant (e.g. its mean) would be suboptimal in separating only alternatives of weighted $L^2$-distance of order $n^{-\alpha/(2\alpha+1)}$. A similar efficiency gain for nonparametric testing is known for Gaussian white noise models, see \cite{Ingster}.

We compare next our test with related existing work.
First, there is an enormous amount of literature on parameter shifts and breaks \citep[and references therein]{Perron2013}, but the results are all based on a long span ergodic-type theory rather than fixed length in-fill conducted here. Second, \cite{AK2012} propose a test for constant beta based on a Hausman type statistic that compares a nonparametric kernel-based estimate of betas at fixed time points and a long-run estimate of beta. \cite{AK2012}  do not consider formally the role of the discretization error in their analysis. By contrast, we rely here solely on a fixed span and the associated in-fill or high-frequency asymptotics, and we are interested in checking whether beta is constant on the whole time interval, not only at fixed points in time. Thus, intuitively, our test checks for constancy of beta on an asymptotically increasing number of blocks of shrinking time span. Third, \cite{TB}, \cite{Kalnina} and \cite{AHHH2012} consider tests for constant integrated betas, i.e., deciding whether integrals of betas over fixed intervals of time such as days or weeks are the same. Unlike these papers, we are interested in deciding whether the spot beta process remains constant within a fixed interval of time which is a stronger hypothesis and requires essentially conducting testing on blocks of shrinking time span. Finally, our work is related to \cite{MZ06}. In the pure diffusive setting (i.e., without jumps), \cite{MZ06} are interested in estimating the residual component of the asset  without any assumption regarding whether the beta remains constant or not while our interest here is in testing the latter.

Finally, our setup is based on equidistant observation grid for the pair $(X,Y)$ and rules out microstructure noise. At ultra-high frequencies asynchronicity and irregularity of sampling times as well as microstructure noise become very important. We believe that our approach can be generalized to accommodate the above features, but the precise technical details will be challenging, see e.g., \cite{Hayashi20112416}, \cite{Bibinger}, \cite{FR} and \cite{bhmr} for the related problem of integrated (multivariate) volatility estimation.

We find satisfactory performance of our estimator on simulated data. In an empirical application we study the appropriate time window width over which market betas of four different assets remain constant. For most of the assets we study we find such a window to be at least as long as a week while for one of the assets our test rejects in a nontrivial number of weeks the null of constancy.

The rest of the paper is organized as follows. In Section~\ref{sec:setup} we introduce our formal setup. In Section~\ref{sec:theory} we develop the test, analyze its behavior under the null and alternative hypothesis, and study its optimality. Section~\ref{sec:mc} contains  a Monte Carlo analysis of the finite sample performance of the test and in Section~\ref{sec:applic} we apply the test to study time-variation of market betas. Section~\ref{sec:concl} concludes. Proofs are in Section~\ref{sec:proofs}. 
\section{Setup and notation}\label{sec:setup}

We start with introducing the setting and stating the assumptions that we need for the results in the paper. The finite variation components of the  underlying bivariate process $(X,Y)$, given in \refeq{eq:XY}, are assumed to be of the form
\begin{equation}\label{eq:XY-2}
\begin{split}
F_t^X = \int_0^t\alpha_s^{X}ds+\int_0^t\int_{E}\delta^X(s,x)\mu(ds,dx),~~F_t^Y =\int_0^t\alpha_s^{Y}ds+\int_0^t\int_{E}\delta^Y(s,x)\mu(ds,dx),
\end{split}
\end{equation}
where $\alpha^{X}$ and $\alpha^{Y}$ are processes with \cadlag paths; $\mu$ is Poisson measure on $\mathbb{R}_+\times \mathbb{R}$ with compensator $dt\otimes dx$; $\delta^X(t,x)$ and $\delta^Y(t,x)$ are two predictable functions. The first and second components of $F^X$ and $F^Y$ in \refeq{eq:XY-2} are the continuous and discontinuous finite variation parts of $X$ and $Y$. We note that for the last integrals in \refeq{eq:XY-2} to make sense, we need jumps to be absolutely summable on finite time intervals. Our setup in \refeq{eq:XY}-\refeq{eq:XY-2}, therefore, implicitly rules out jumps of infinite variation. This is similar to prior work on estimation of integrated volatility because infinite variation jumps necessarily spoil inference on the diffusion part of the processes, cf. \cite{jacodreiss}.

We further assume that the volatility processes $\sigma$ and $\widetilde{\sigma}$ are themselves \Ito semimartingales, i.e., they have representations of the form
\begin{equation}\label{eq:vol}
\begin{split}
\sigma_t&=\sigma_0+\int_0^t\alpha_s^{\sigma}ds+\int_0^t\gamma_s^{\sigma}dW_s+\int_0^t\widetilde{\gamma}_s^{\sigma}d\widetilde{W}_s+\int_0^t\gamma_s^{'}dW_s^{'}
\\&~~~+\int_0^t\int_E\kappa(\delta^{\sigma}(s,x))\widetilde{\mu}(ds,dx)+\int_0^t\int_E\kappa'(\delta^{\sigma}(s,x))\mu(ds,dx),\\
\widetilde{\sigma}_t&=\widetilde{\sigma}_0+\int_0^t\alpha_s^{\widetilde{\sigma}}ds+\int_0^t\gamma_s^{\widetilde{\sigma}}dW_s+\int_0^t
\widetilde{\gamma}_s^{\widetilde{\sigma}}d\widetilde{W}_s+\int_0^t\gamma_s^{''}dW_s^{''}\\&~~~+\int_0^t\int_E\kappa(\delta^{\widetilde{\sigma}}(s,x))\widetilde{\mu}(ds,dx)
+\int_0^t\int_E\kappa'(\delta^{\widetilde{\sigma}}(s,x))\mu(ds,dx),
\end{split}
\end{equation}
where $W'$ and $W^{''}$ are two Brownian motions, having arbitrary dependence, but independent from $(W_t,\widetilde{W}_t)$; $\widetilde{\mu}(dt,dx) = \mu(dt,dx)-dt\otimes dx$ is the compensated jump measure; $\kappa(\cdot)$ is a continuous function with bounded domain and with $\kappa(x)=x$ in a neighborhood of zero, $\kappa'(x)=x-\kappa(x)$; $\alpha^{\sigma}$, $\alpha^{\widetilde{\sigma}}$, $\gamma^{\sigma}$, $\gamma^{\widetilde{\sigma}}$, $\widetilde{\gamma}^{\sigma}$, $\widetilde{\gamma}^{\widetilde{\sigma}}$, $\gamma^{'}$ and $\gamma^{''}$ are processes with \cadlag paths; $\delta^{\sigma}(t,x)$ and $\delta^{\widetilde{\sigma}}(t,x)$ are two predictable functions.

We note that the specification in \refeq{eq:XY}-\refeq{eq:vol} is very flexible and allows for most of the stochastic volatility models considered in empirical work. We also allow for arbitrary dependence between the Brownian motion and Poisson measure driving $X$ and the volatility processes. We state our assumptions for \refeq{eq:XY}-\refeq{eq:vol} in the following.

\noindent \textbf{Assumption A.} For the process defined in \refeq{eq:XY}-\refeq{eq:vol} we have:
\begin{itemize}
\item [(a)] \textit{$|\sigma_t|^{-1}$, $|\sigma_{t-}|^{-1}$, $|\widetilde{\sigma}_t|^{-1}$ and $|\widetilde{\sigma}_{t-}|^{-1}$ are strictly positive; }
\item [(b)] \textit{$\beta$, $\alpha^{\sigma}$, $\alpha^{\widetilde{\sigma}}$, $\gamma^{\sigma}$, $\gamma^{\widetilde{\sigma}}$, $\widetilde{\gamma}^{\sigma}$, $\widetilde{\gamma}^{\widetilde{\sigma}}$, $\gamma'$ and $\gamma^{''}$ are \cadlag adapted; $\delta^X$, $\delta^Y$, $\delta^{\sigma}$ and $\delta^{\widetilde{\sigma}}$ are predictable;}
\item [(c)] \textit{$\alpha^X$ and $\alpha^Y$ are \Ito semimartingales with locally bounded coefficients;}
\item [(d)] \textit{There is a sequence $T_k$ of stopping times increasing to infinity such that:}
\begin{equation}
t\leq T_k~~\Longrightarrow~~|\delta^{X}(t,x)|\wedge 1+|\delta^{Y}(t,x)|\wedge 1\leq \gamma_k^{(1)}(x),~~|\delta^{\sigma}(t,x)|\wedge 1+|\delta^{\widetilde{\sigma}}(t,x)|\wedge 1\leq \gamma_k^{(2)}(x),\nonumber
\end{equation}
\textit{where $\gamma_k^{(1)}(x)$ and $\gamma_k^{(2)}(x)$ are deterministic functions on $\mathbb{R}$ satisfying}
\begin{equation}
\int_{\mathbb{R}}|\gamma_k^{(1)}(x)|^rdx<\infty,~~\textrm{and}~~\int_{\mathbb{R}}|\gamma_k^{(2)}(x)|^2dx<\infty,\nonumber
\end{equation}
\textit{for some $r\in[0,1]$.}
\end{itemize}

Parts (a) and (b) of Assumption A are necessary as our inference on $\beta_t$ depends on the presence of the diffusion components in $X$ and $Y$. Part (c) of Assumption A controls the activity of the jumps in $X$ and $Y$ and some of our results will depend on the number $r$.


\section{Main results}\label{sec:theory}

We proceed with formulating the testing problem that we study in the paper. We assume that we observe the process $(X,Y)$ on the interval $[0,1]$ at the equidistant grid $0,\frac{1}{n},\frac{2}{n},...,1$ for some $n\in\mathbb{N}$, and the asymptotics in the paper will be for $n\rightarrow\infty$. The results, of course, extend trivially to arbitrary time intervals of fixed length. Our interest lies in designing a test to decide whether the stochastic spot beta process $\beta$ remains constant or not on the interval $[0,1]$. This is a pathwise property and therefore we are interested in discriminating the following two events dividing the sample space:
\begin{equation}
\Omega^c = \left\{\omega: \beta_t(\omega)=\beta_0(\omega)~~\textrm{almost everywhere on}~ [0,1]   \right\},~~~\Omega^v = \Omega\setminus\Omega^c.
\end{equation}
The set $\Omega^c$ can be characterized in different ways. One natural way is
\begin{equation}\label{eqOmegac}
\Omega^c = \left\{\omega: \int_0^1\beta_s^2(\omega)\,ds-\left(\int_0^1\beta_s(\omega)\,ds\right)^2=0 \right\},
\end{equation}
which in words means that $\beta_t$ is constant on the interval $[0,1]$ if and only if its variance on that interval with respect to the occupation measure associated with $\beta$ vanishes. One can then formulate a feasible test by constructing estimates for $\int_0^1\beta_s^2ds-\left(\int_0^1\beta_sds\right)^2$ from the high-frequency data on $(X,Y)$. This can be done by forming blocks with increasing number of observations in each of them but with shrinking time span and estimating $\beta_t$ locally in each of the blocks, following a general approach proposed in \cite{JR} (see also \cite{MZ06}). It turns out, however, that under the null of constant beta, a CLT for $\int_0^1\beta_s^2ds-\left(\int_0^1\beta_sds\right)^2$ as in \cite{JR} is degenerate and higher order asymptotics is needed. This is because the derivatives of the test statistic with respect to the elements of the variance-covariance matrix on the blocks, used to construct an estimate for $\int_0^1\beta_s^2ds-\left(\int_0^1\beta_sds\right)^2$, vanish under the null hypothesis. Besides, in this case we also need  debiasing terms.

Therefore, we adopt here an alternative point of view to characterize $\Omega^c$ that avoids the above complications. Suppose that we know the value of $\beta_t$ at time $t=0$. In this case, recalling that the process $\sigma$ is non-vanishing on the interval $[0,1]$, we have
\begin{equation}
\Omega^c = \left\{\omega: \langle Y^c-\beta_0X^c,X^c\rangle_t = 0, ~~\textrm{for every $t\in[0,1]$}\right\}.
\end{equation}
If $X^c$ and $Y^c$ had constant and deterministic volatility, we would have to test for independence in the bivariate Gaussian sample $(\Delta_i^nX^c,\Delta_i^n(Y^c-\beta_0X^c))_{1\le i\le n}$, where henceforth we use the shorthand $\Delta_i^nZ = Z_{\frac{i}{n}}-Z_{\frac{i-1}{n}}$ for an arbitrary process $Z$. In this case, the natural (i.e., uniformly most powerful unbiased) test is of the form $nR^2-1>c$ with the sample correlation coefficient
\begin{equation}\label{eq:R}
R=\frac{\sum_i \Delta_i^nX^c\Delta_i^n(Y^c-\beta_0X^c)}{\sqrt{\sum_i (\Delta_i^n X^c)^2}\sqrt{\sum_i (\Delta_i^n(Y^c-\beta_0X^c))^2}}.
\end{equation}
The critical value $c>0$ is distribution-free and derived from the finite sample result that $\sqrt{n-1}R/\sqrt{1-R^2}$ follows a $t_{n-1}$-distribution under the independence hypothesis (this follows from the exact finite sample distribution of regressions with normal errors, see e.g., chapter 1 of \cite{Hayashi}) or asymptotically from $nR^2\to \chi^2(1)$. Since in our case, the volatilities are time varying we base our test on localised statistics of this kind. In line with optimal testing  for nonparametric regression functions, cf. Section \ref{SecNonparTest} below, the final test is based on the sum of these localised test statistics which guarantees high power against time-varying  $\beta_t$ deviating from $\beta_0$ in terms of a weighted $L^2$-distance.  With this in mind, we turn to the concrete construction of our test statistics.

\subsection{The test statistics}

We split the high-frequency observations into $\lfloor n/k_n\rfloor$ blocks with $k_n$ observations per block for $k_n\rightarrow\infty$ and $k_n/n\rightarrow0$. For some constant $b\in\mathbb{R}$, $\alpha>0$ and $\varpi\in(0,1/2)$, we introduce
\begin{equation}
\widehat{C}_j^n(b) = \frac{n}{\sqrt{k_n}}\sum_{i=(j-1)k_n+1}^{jk_n}\Delta_i^nX(\Delta_i^nY-b\Delta_i^nX)1_{\{|\Delta_i^nX|\leq\alpha\Delta_n^{\varpi},~|\Delta_i^nY|\leq\alpha
\Delta_n^{\varpi}\}}.
\end{equation}
$\widehat{C}_j^n(b)/\sqrt{k_n}$ is an estimate of $\frac{n}{k_n}\left[\langle Y^c-b X^c,X^c\rangle_{\frac{jk_n}{n}}-\langle Y^c-b X^c,X^c\rangle_{\frac{(j-1)k_n}{n}}\right]$ which is zero if the continuous beta is constant and $b=\beta_0$. In this case, because of the shrinking time span of the block, for our purposes $\widehat{C}_j^n(\beta_0)$ will be equivalent to $\sigma_{\frac{(j-1)k_n}{n}}\widetilde{\sigma}_{\frac{(j-1)k_n}{n}}\frac{n}{\sqrt{k_n}}\sum_{i=(j-1)k_n+1}^{jk_n}\Delta_i^nW\Delta_i^n\widetilde{W}$ asymptotically. So, conditionally on ${\cal F}_{\frac{(j-1)k_n}{n}}$ we are in the above bivariate Gaussian setting.

The analogue of the denominator of $R^2$ in \eqref{eq:R} is given by
\begin{equation}
\widehat{V}_j^n(b) = \widehat{V}_j^{(n,1)}\widehat{V}_j^{(n,2)}(b),~~\widehat{V}_j^{(n,1)}= \frac{n}{k_n}\sum_{i =
(j-1)k_n+1}^{jk_n}\left(\Delta_i^nX\right)^21_{\{|\Delta_i^nX|\leq\alpha\Delta_n^{\varpi},~|\Delta_i^nY|\leq\alpha\Delta_n^{\varpi}\}},
\end{equation}
\begin{equation}
\widehat{V}_j^{(n,2)}(b) = \frac{n}{k_n}\sum_{i=(j-1)k_n+1}^{jk_n}\left(\Delta_i^nY-b\Delta_i^nX\right)^21_{\{|\Delta_i^nX|\leq\alpha\Delta_n^{\varpi},~|\Delta_i^nY|\leq\alpha
\Delta_n^{\varpi}\}}.
\end{equation}
Here, however, we compensate $\widehat C_j^n(b)^2$ by $\widehat{V}_j^n(b)$ in the numerator, while dividing by
the estimate from the previous block, $\widehat{V}_{j-1}^n(b)$. The predictable choice of the denominator (with respect to ${\cal F}_{\frac{(j-1)k_n}{n}}$) guarantees a block-wise martingale difference property and thus avoids an additional bias in the case of stochastic volatilities.  Thus, our test statistic takes the final form
\begin{equation}
\widehat{T}^n(b)
= \frac{1}{\sqrt{2}}\sqrt{\frac{k_n}{n}}\sum_{j=2}^{\lfloor\frac{n}{k_n}\rfloor}\widehat{T}_j^n(b),~~\widehat{T}_j^n(b) =  \frac{\left(\widehat{C}_j^n(b)\right)^2-\widehat{V}_j^n(b)}{\widehat{V}_{j-1}^n(b)}.
\end{equation}
Let us point out that for convenience all statistics $\widehat{C}_j^n(b), \widehat{V}_j^n(b), \widehat{T}_j^n(b)$ and $\widehat{T}^n(b)$ are scaled to have stochastic order one under the null of constant beta.

\subsection{Testing for a known constant beta}

We proceed next with deriving the limit behavior of our statistic in the case when testing for a known constant beta. Our test statistic $\widehat{T}^n(b)$ is a sum of nonlinear transforms of block-based volatility estimates. \cite{JR} derive the limit behavior of statistics of this type (for some sufficiently smooth functions of the volatility estimates), with the asymptotic limit being determined by a first-order linear approximation of the function around the volatility level over the (shrinking) block. In our case such a local linearization approach does not work because the block-based covariance estimate $\widehat{C}_j^n(\beta_0)/\sqrt{k_n}$ converges to zero and hence we need to scale it up (by $\sqrt{k_n}$). To derive the limit behavior of $\widehat{T}^n(\beta_0)$ here we directly approximate the conditional moments of $\widehat{T}_j^n(\beta_0)$ and further make use of the fact that $\{\widehat{T}_j^n(\beta_0)\}_{j=2,...,\lfloor\frac{n}{k_n}\rfloor}$ are approximately mean zero and uncorrelated across blocks. The formal result is given in the next theorem.

\begin{thm}\label{thm:known_beta}
Grant Assumption A and let the sequence $(k_n)$ satisfy $k_n\rightarrow \infty$ with $\frac{k_n}{n}\rightarrow 0$.
\begin{enumerate}
\item If $k_n^{-1}n^{1/4}\rightarrow 0$ and $k_n^{-1}n^{2-(4-r)\varpi}\rightarrow 0$ with $\varpi\in\left(\frac{1}{2(2-r)},\frac{1}{2}\right)$,  we have
\begin{equation}\label{eq:T_null_known}
\widehat{T}^n(\beta_0)~\stackrel{\mathcal{L}}{\longrightarrow}~Z,~~~\textrm{in restriction to the set $\Omega^c$,}
\end{equation}
for $Z$ being a standard normal random variable.

\item  If $k_n^{-1}n^{1-(2-r)\varpi}\rightarrow 0$,  we have
\begin{equation}\label{eq:T_alter_known}
\frac{1}{\sqrt{nk_n}}\widehat{T}^n(\beta_0)~\stackrel{\mathbb{P}}{\longrightarrow}~\frac{1}{\sqrt{2}}\int_0^1\frac{(\beta_{s}-\beta_0)^2\sigma_{s}^2}
{\left((\beta_{s}-\beta_0)^2\sigma_{s}^2
+\widetilde{\sigma}^2_{s}\right)}ds,~~~\textrm{in restriction to the set $\Omega^v$.}
\end{equation}
\end{enumerate}
\end{thm}
Starting with the behavior under the null hypothesis of constant beta, we see that the asymptotic limit of our statistic is standard normal and does not depend on any of the ``nuisance parameters'' in our model like the volatility processes $\sigma$ and $\widetilde{\sigma}$. This is due to the fact that the statistic is of ``self-normalizing'' type. This is very convenient for the inference process. In addition, the self-normalization property of our statistic avoids the need of showing stable convergence (which is a much stronger form of convergence), typically needed in high-frequency asymptotics for conducting feasible inference, see e.g., \cite{JP12}. The condition on the block size $k_n^{-1}n^{1/4}\rightarrow 0$ in part(a) of Theorem~\ref{thm:known_beta} is to ensure that the averaging within the blocks is sufficient so that the within-block averages are not far away from their limits. The condition $k_n^{-1}n^{2-(4-r)\varpi}\rightarrow 0$ is to ensure that the error due to the elimination of the jumps is negligible. The user chooses $\varpi$, so as with estimators of truncated type (\cite{Ma1}), it is optimal to set $\varpi$ as close as possible to its upper limit of $1/2$. In this case, the lower bound on $\varpi$ in Theorem~\ref{thm:known_beta}(a) will be satisfied (provided $r<1$). We note also that the second condition for $k_n$ in part(a) of the theorem becomes more restrictive  for higher values of the jump activity as the separation of higher activity jumps from the diffusive component is harder.

Turning to the limit of our statistic in the case of $\beta$ time-varying on the interval $[0,1]$, given in part (b) of the theorem, we see that the limit is a weighted average of the distance $(\beta_s-\beta_0)^2$. The weighting is determined by the stochastic volatilities $\sigma_s^2$ and $\widetilde{\sigma}_s^2$ over the interval. The scaling down of the statistic is by the factor $\sqrt{nk_n}$, which means that higher block size $k_n$ leads to higher rate of explosion of the statistic under the alternative. Finally, the condition for the block size in part(b) of the theorem is very close to the analogous one under the null hypothesis in part (a) of the theorem, provided $\varpi$ is selected very close to $1/2$.

The limit in \refeq{eq:T_alter_known} reveals the difficulties in detecting time variation in beta with our test. In particular, keeping everything else fixed, a higher level of the idiosyncratic volatility, $\widetilde{\sigma}^2$, decreases the value of the limit in \refeq{eq:T_alter_known}, and hence reduces our ability to detect time variation in beta. The effect of the systematic volatility, $\sigma^2$, is in the opposite direction. A higher value of the systematic volatility means that the systematic component of $Y$ has a bigger share in its total variation.

\subsection{Testing for unknown constant beta}

In most cases of practical interest, we will not know the level of beta, but instead we shall need to estimate it under the assumption that it is constant over a given interval. We will then be simultaneously interested in the estimated value and in the outcome of a test to decide whether it can be assumed to have stayed constant. Thus, we need first an initial estimator of the continuous beta over the interval. We shall use the following natural estimator
\begin{equation}\label{beta_hat}
\widehat{\beta}_n = \frac{\sum_{i=1}^{n}\Delta_i^nX\Delta_i^nY1_{\{|\Delta_i^nX|\leq\alpha\Delta_n^{\varpi},~|\Delta_i^nY|\leq\alpha\Delta_n^{\varpi}\}}}
{\sum_{i=1}^{n}(\Delta_i^nX)^21_{\{|\Delta_i^nX|\leq\alpha\Delta_n^{\varpi},~|\Delta_i^nY|\leq\alpha\Delta_n^{\varpi}\}}},
\end{equation}
which can be equivalently defined as
\begin{equation}
\widehat{\beta}_n = \textrm{argmin}_{\beta}\sum_{i=1}^{n}\left(\Delta_i^nY-\beta\Delta_i^nX\right)^21_{\{|\Delta_i^nX|\leq\alpha\Delta_n^{\varpi},~|\Delta_i^nY|\leq\alpha\Delta_n^{\varpi}\}},
\end{equation}
where the objective function in the above optimization is the empirical analogue of $\langle Y^c-\beta X^c, Y^c-\beta X^c\rangle_1$. This estimator has been studied in \cite{TB} and \cite{GM}.

When the process $\beta$ varies over the time interval $[0,1]$, $\widehat{\beta}_n$ converges in probability to
\begin{equation}
\overline{\beta} = \frac{ \int_0^1\beta_s\sigma_s^2ds } {\int_0^1\sigma_s^2ds},
\end{equation}
which can be viewed as a volatility weighted average of the time-varying beta over the interval. The rate of convergence of $\widehat{\beta}_n$ is $\sqrt{n}$ and its limiting behavior in the general case when the process $\beta$ can vary over time is given by the following lemma.

\begin{lema}\label{lema:beta_hat}
Suppose the process $(X,Y)$ satisfies Assumption A and let $\varpi\in\left(\frac{1}{2(2-r)},\frac{1}{2}\right)$. Then
\begin{equation}
\sqrt{n}\left(\widehat{\beta}_n - \overline{\beta}\right)~\stackrel{\mathcal{L}-s}{\longrightarrow}~\sqrt{V_{\beta}}Z,
\end{equation}
where $Z$ is independent standard normal defined on an extension of the original probability space and
\begin{equation}
\begin{split}
V_{\beta} &= \frac{2}{\left(\int_0^1\sigma_s^2ds\right)^4}\bigg[\int_0^1(\beta_s^2\sigma_s^4+\sigma_s^2\widetilde{\sigma}_s^2)ds\left(\int_0^1\sigma_s^2ds\right)^2+
\left(\int_0^1\beta_s\sigma_s^2ds\right)^2\int_0^1\sigma_s^4ds
\\&~~~~~~~~~~~~~~~~~~~~~~~~~~~~~-2\int_0^1\beta_s\sigma_s^2ds\int_0^1\beta_s\sigma_s^4ds\int_0^1\sigma_s^2ds\bigg].
\end{split}
\end{equation}
\end{lema}
The proof of Lemma~\ref{lema:beta_hat} follows from the limiting results for multivariate truncated variation, see e.g., Theorem 13.2.1 of \cite{JP12},  and an application of the Delta method.

With this estimator of $\beta_0$ (under the null), our test in the case of unknown beta is simply based on $\widehat{T}^n(\widehat{\beta}_n)$. Its asymptotic behavior is given in the following theorem.

\begin{thm}\label{thm:unknown_beta}
Grant Assumption A and let the sequence $(k_n)$ satisfy $k_n\rightarrow \infty$ with $\frac{k_n}{n}\rightarrow 0$.
\begin{enumerate}
\item  If $k_n^{-1}n^{1/4}\rightarrow 0$ and $k_n^{-1}n^{2-(4-r)\varpi}\rightarrow 0$ with $\varpi\in\left(\frac{1}{2(2-r)},\frac{1}{2}\right)$,  we have
\begin{equation}\label{eq:T_null_unknown}
\widehat{T}^n(\widehat{\beta}_n)~\stackrel{\mathcal{L}}{\longrightarrow}~Z,~~\textrm{in restriction to the set $\Omega^c$,}
\end{equation}
for $Z$ being a standard normal random variable.

\item If $k_n^{-1}n^{1-(2-r)\varpi}\rightarrow 0$,  we have
\begin{equation}\label{eq:T_alter_unknown}
\frac{1}{\sqrt{nk_n}}\widehat{T}^n(\widehat{\beta}_n)~\stackrel{\mathbb{P}}{\longrightarrow}~\frac{1}{\sqrt{2}}\int_0^1\frac{(\beta_{s}-\overline{\beta})^2\sigma_{s}^2}
{\left((\beta_{s}-\overline{\beta})^2
\sigma_{s}^2+\widetilde{\sigma}^2_{s}\right)}ds,~~\textrm{in restriction to the set $\Omega^v$.}
\end{equation}
\end{enumerate}
\end{thm}

From part(a) of the theorem we can see that the estimation of the unknown beta has no asymptotic effect on our statistic under the null. The only difference from the testing against a known constant beta under the alternative is that now the limit of the statistic in \refeq{eq:T_alter_unknown} contains the averaged value $\overline{\beta}$. Note that the limit of \refeq{eq:T_alter_unknown} is a volatility weighted version of \refeq{eqOmegac}.

\subsection{Testing against local alternatives and asymptotic optimality}\label{SecNonparTest}

The asymptotics under the alternative in Theorems \ref{thm:known_beta}(b) and \ref{thm:unknown_beta}(b) are somewhat misleading regarding the choice of the block size $k_n$. For a fixed single alternative the test is asymptotically most powerful if $k_n$ is chosen as large as possible. This, however, is not reasonable for fixed $n$ because on large blocks time varying betas that oscillate will give similar values for the test statistics as constant betas due to the averaging on each block. This phenomenon is well understood for testing a nonparametric regression function where the bandwidth $h$ of a kernel smoother takes on the role of the relative block size $k_n/n$. For a more meaningful statement local alternatives as well as uniform error probabilities should be considered.

Following \cite{Ingster} we are studying the optimal separation rate $r_n$ between the single hypothesis
\[H_0=\{P_\beta\},\]
for some fixed constant risk value $\beta>0$ and the local nonparametric alternative
\[H_{1,\alpha}(r_n)=\Big\{P_{\beta_t}\text{ such that a.s. }\beta_t\in C^\alpha(R), \, \int_0^1\frac{{\sigma_t^2(\beta_t-\beta)^2}}{\sigma_t^2(\beta_t-\beta)^2+\widetilde{\sigma}_t^2}\,dt\ge r_n^2,\Big\},
\]
where $C^\alpha(R)=\{f\,:\,|f(t)-f(s)|\le R|t-s|^\alpha,\,|f(t)|\le R\}$ for all $t,s\in[0,1]$ denotes a H\"older ball of regularity $\alpha\in(0,1]$ and radius $R>0$. In this notation it is understood that the laws $P_{\beta_t}$ are defined on the path space of $((X_t,Y_t),t\in[0,1])$ and the nuisance parameters $\sigma_t^2,\widetilde{\sigma}_t^2$ and the drift and jump parts may vary with the parameter of interest $\beta_t$.

The separation rate $r_n\downarrow 0$ is called minimax optimal over $C^\alpha(R)$ if there is a  test $\phi_n$, based on $n$ observations, such that
\[ \forall \gamma\in(0,1)~~\exists \Gamma>0:\;\limsup_{n\to\infty}\Big(P_\beta(\phi_n=1)+\sup_{P_{\beta_t}\in H_{1,\alpha}(\Gamma r_n)}P_{\beta_t}(\phi_n=0)\Big)\le\gamma\]
holds while the infimum of the error probabilities over any possible test $\psi_n$ remains positive:
\[ \forall \gamma\in(0,1)~~\exists \tilde\Gamma>0:\liminf_{n\to\infty}\inf_{\psi_n}\Big(P_\beta(\psi_n=0)+\sup_{P_{\beta_t}\in H_{1,\alpha}(\tilde\Gamma r_n)}P_{\beta_t}(\psi_n=0)\Big)\ge\gamma.
\]
Our test then satisfies a minimax bound with separation rate $r_n=n^{-2\alpha/(4\alpha+1)}$. To keep the proofs transparent, we show this only in the case when $X$ and $Y$ do not jump and $\beta$ from the hypothesis $H_0$ is assumed to be known.

\begin{thm}\label{thm:la_sup}
Assume that Assumption SA in Section~\ref{subsec:proofs_loc} holds and $\delta^X(t,x) = \delta^Y(t,x)=0$ for $t\in[0,1]$. Suppose $\alpha>5/12$ and $k_n=\lfloor n^{\frac{4\alpha-1}{4\alpha+1}}\rfloor$, $r_n=n^{-2\alpha/(4\alpha+1)}$. Then for any $\gamma\in(0,1)$ and critical value $c_{\gamma/2}$ under the hypothesis (i.e. $\limsup_{n\to\infty}P_{\beta}(\widehat T^{n}(\beta)\ge c_{\gamma/2})\le\gamma/2$), there is a $\Gamma>0$ such that the test $\phi_n=1\{\widehat T_n(\beta)>c_{\gamma/2}\}$ satisfies
\begin{equation}
\limsup_{n\to\infty}\Big(P_\beta(\phi_n=1)+\sup_{P_{\beta_t}\in H_{1,\alpha}(\Gamma r_n)}P_{\beta_t}(\phi_n=0)\Big)\le \gamma.
\end{equation}
\end{thm}
The condition $\alpha>5/12$ in the above theorem  is due to the rate condition on the block size $k_n^{-1}n^{1/4}\to 0$ in Theorem~\ref{thm:known_beta}(a).
The natural assumption for the process $\beta$ is that it is itself a continuous It\^o semimartingale and thus has H\"older regularity $\alpha$ of almost $1/2$. In this case the optimal block length is $k_n\approx n^{1/3}$ and the separation rate is $r_n\approx n^{-1/3}$, which is far better than the optimal nonparametric estimation rate $n^{-\alpha/(2\alpha+1)}\approx n^{-1/4}$.

For an adaptive choice of the optimal block length $k_n$, without specifying the regularity $\alpha$ in advance, and even ``parametric power'' (in the sense of Theorem~\ref{thm:known_beta}(b)) for certain parametric submodels for $\beta_t$, an analogue of the maximal test statistics of \cite{horowitz2001adaptive}  can be applied. Note that they also show that their test allows for a parametric form of the null hypothesis, assuming that the true parameter can be estimated at rate $n^{-1/2}$ under the null. This estimator is plugged into the test statistics exactly in the same way as we test for unknown $\beta$.

Here, we focus on the non-obvious question of optimality. We shall derive a lower bound on the separation rate for the even smaller subclass of pure Gaussian martingales, where neither jumps nor drift terms appear in $(X,Y)$ and where the volatilities are deterministic.
Already in this subclass no other test can have a smaller minimax separation rate than $r_n=n^{-2\alpha/(4\alpha+1)}$, which then, of course, extends to the more general model for which our test is designed. Our test is thus indeed minimax optimal.
\begin{thm}\label{thm:la_inf}
Assume that $\sigma_t^2$ and $\widetilde{\sigma}_t^2$ are $\alpha$-\Holder~continuous deterministic functions, for some $\alpha\in(0,1]$, and $\alpha^X_t=\alpha^Y_t=\delta^X(t,x)=\delta^Y(t,x)=0$ for $t\in[0,1]$. Then
for any $\alpha\in(0,1]$, $\gamma\in(0,1)$ there is a $\tilde\Gamma>0$ such that for $r_n=n^{-2\alpha/(4\alpha+1)}$ and arbitrary tests $\psi_n$
\begin{equation}
\liminf_{n\to\infty}\inf_{\psi_n}\Big(P_\beta(\psi_n=1)+\sup_{P_{\beta_t}\in H_{1,\alpha}(\tilde\Gamma r_n)}P_{\beta_t}(\psi_n=0)\Big)\ge\gamma.
\end{equation}
\end{thm}

\section{Monte Carlo study}\label{sec:mc}

We now evaluate the performance of our test on simulated data from the following model
\begin{equation}\label{eq:mc_model}
\begin{split}
dX_t &= \sqrt{V_t}dW_t+dL_t,~~dY_t = \beta_tdX_t+\sqrt{\widetilde{V}_t}d\widetilde{W}_t+d\widetilde{L}_t,\\
dV_t &= 0.03(1-V_t)dt+0.18\sqrt{V_t}dB_t,~~d\widetilde{V}_t = 0.03(1-\widetilde{V}_t)dt+0.18\sqrt{\widetilde{V}_t}d\widetilde{B}_t,\\
\end{split}
\end{equation}
where $(W,\widetilde{W},B,\widetilde{B})$ is a vector of independent standard Brownian motions; $L$ and $\widetilde{L}$ are two  pure-jump \Lvb processes, independent of each other and of the Brownian motions, each of which with characteristic triplet $(0,0,\nu)$ for a zero truncation function and $\nu(dx) = 1.6e^{-2|x|}dx$. $V$ and $\widetilde{V}$ in \refeq{eq:mc_model} are square-root diffusion processes used extensively in financial applications for modeling volatility. For the process $\beta$, we consider
\begin{equation}
H_0:\beta_t=1~~\textrm{and}~~H_a:~~d\beta_t = 0.03(1-\beta_t)dt+0.18\sqrt{\beta_t}dB_t^{\beta},
\end{equation}
for $B^{\beta}$ being a Brownian motion independent from the Brownian motions in \refeq{eq:mc_model}. The parameters of the model are calibrated to real financial data. In particular, the means of $V_t$ and $\widetilde{V}_t$ are set to $1$ and they are both persistent processes (our unit of time is a trading day and returns are in percentage). Jumps in $X$ and $Y$ have intensity of $0.4$ jumps per day and $0.8$ jumps per day respectively. The variances of the jump components of both $X$ and $Y$ are $40\%$ that of their continuous components (on any fixed time interval).

The observation scheme is similar to that of our empirical application. We set $1/\Delta_n = 38$, which corresponds to sampling every $10$ minutes in a $6.5$ hours trading day. In the application of the test, we set $k_n=19$ which corresponds to constructing two blocks per unit of time (which is day). We test for constant beta on an interval of length of $T=5$ (week), $T=22$ (one month) and $T=66$ (one quarter) by summing the test statistics over the $T$ days.

The results from the Monte Carlo, which is based on $1000$ replications,  are reported in Table~\ref{tb:mc}. The test performs reasonably in finite samples. In particular, the actual rejection rates are in the vicinity of the nominal ones under the null hypothesis of constant beta across the three intervals $T=5$, $T=22$ and $T=66$. We notice a bit of over-rejection at the $1\%$ level across the three intervals. Turning to the power of the test, not surprisingly we note that the power increases with $T$, with the power against the considered time-varying beta model being lowest for the case $T=5$. Intuitively, more observations (higher $T$) allow us to better discriminate the noise in the recovery of $\beta$ from its true time variation.

\begin{table}[h]
\begin{center}
\begin{tabular*}{4.2in}{crrrrrrr} \hline\hline
\textbf{Interval} &  \multicolumn{7}{c}{\textbf{ Significance Level (Percent)}}\\\\[-0.25ex]
                  &  $10.0$ & $5.0$ & $1.0$ & & $10.0$ & $5.0$ & $1.0$\\[+1.25ex]
&  \multicolumn{3}{c}{\textbf{Constant Beta }} & & \multicolumn{3}{c}{\textbf{Time-Varying Beta }} \\[+0.75ex]
week    &  $7.11$  & $4.70$  & $2.30$ &  & $12.64$ & $8.95$  & $5.07$\\
month   & $10.50$  & $6.30$  & $3.00$ &  & $45.97$ & $39.43$ & $29.17$\\
quarter & $10.70$  & $7.10$  & $3.10$ &  & $83.20$ & $79.50$ & $72.20$\\
 \hline\hline
\end{tabular*}
\end{center}
\caption{Monte Carlo Results.}
\label{tb:mc}
\end{table}

\section{Empirical application}\label{sec:applic}

The test for constant market beta is conducted on four assets sampled at the 10-minute frequency over the period $2006$--$2012$.  We refer to them by ticker symbol: IBM, XOM (Exxon Mobil), GLD (an Exchange-traded Fund (ETF) that tracks the price of gold), and BAC (Bank of America). IBM and XOM are both very stable large-cap stocks; GLD (or gold) is a storable asset that provides a hedge against general macroeconomic risks, while BAC went through stressful episodes with large price fluctuations during the global financial crises.  The market index is $SPY$, the ETF that tracks the S\&P 500 index.

Each 10-minute data set consists of $1746$ days of $38$ within-day returns (log-price increments), and the tests are conducted at the weekly, monthly, and quarterly time intervals.  A week consists of five consecutive trading days, while the calendar months and quarters contain (on average) $22$ and $66$ trading days, respectively. We use the term windows for these time segments. The test is implemented exactly as in the Monte Carlo, in particular we set the block size to $k_n=19$.

Table~\ref{tb:weekly_beta} shows the observed rejection rates of the test for constant, but unknown, beta over the three windows for different size levels and each of the four securities. Starting with IBM, for the weekly window there is little evidence against the null of constant beta at the 10 and 5 percent levels and only slightly so at the 1 percent level (but recall from the Monte Carlo that at $1\%$ our test is slightly over-rejecting in finite samples). On the other hand, the observed rejection rates are somewhat above nominal for a monthly window and well above nominal for a quarter interval. We detect a very similar pattern for XOM. Mainly, at the weekly window there is no strong statistical evidence for time-varying betas while the rejection rates of the test for constant betas increases well above nominal levels as we move from a monthly to a quarterly window. We note that the evidence for time-variation in the market beta of XOM at the monthly and quarterly level is quite stronger than that for IBM. Interestingly for GLD the results are much the same, despite the fact that gold is just a storable commodity with negative cost of carry and used largely as a reserve asset in contrast to IBM and XOM, two huge profitable enterprizes. Taken together, the results suggest that for IBM, XOM and GLD, a weekly window would be a safe choice for treating market beta as constant in an asset pricing study.

On the other hand, the conclusions from Table~\ref{tb:weekly_beta} for BAC  are far different. The betas appear unstable for any testing window at all three nominal frequencies.  In retrospect, this instability might not be surprising given the changing corporate structure and regulatory environment experienced by this company over the period $2006$---$2012$. The outcomes in the table suggest it would be misguided and perhaps misleading to undertake an asset pricing test of BAC treating its market beta as constant over any of the considered windows.

\begin{table}[t]
\begin{center}
\begin{tabular*}{4.2in}{rrrrrrrrr} \hline\hline
\textbf{Interval} &  \multicolumn{6}{c}{\textbf{ Significance Level (Percent)}}\\\\[-2.00ex]
                  &  $10.0$ & $5.0$ & $1.0$ && $10.0$ & $5.0$ & $1.0$ \\[+1.25ex]
 & \multicolumn{3}{c}{\textbf{ IBM }} && \multicolumn{3}{c}{\textbf{ XOM  }} \\[+0.75ex]
week       &     10.60 &   6.59 &   3.72  & &    10.60  &  8.02  &  3.44\\
month      &     20.24 &  10.71 &   7.14  & &    29.76  & 21.43  & 16.67\\
quarter    &     42.86 &  39.29 &  14.29  & &    57.14  & 50.00  & 28.57\\
  \\
 & \multicolumn{3}{c}{\textbf{ GLD }} && \multicolumn{3}{c}{\textbf{ \hspace{-0.10in} BAC}} \\[+0.75ex]
week       &      7.74 &   5.16 &    2.87 & &    14.61   & 10.89 &   6.02\\
month      &     22.62 &  16.67 &    9.52 & &    38.10   & 30.95 &  23.81\\
quarter    &     64.29 &  50.00 &   39.29 & &    78.57   & 71.43 &  57.14\\
\hline\hline
\end{tabular*}
\end{center}
\caption{Tests for Constant Market Betas. \textit{See text for securities associated with the ticker symbols.  For each specified window length, the table shows the percent of all windows for
which the hypothesis of constant but unknown beta is rejected at the specified nominal level.}}\label{tb:weekly_beta}
\end{table}

\begin{figure}[t]
\begin{center}
\begin{tabular}{c}
\includegraphics[width=0.8\textwidth, height=4.5in]{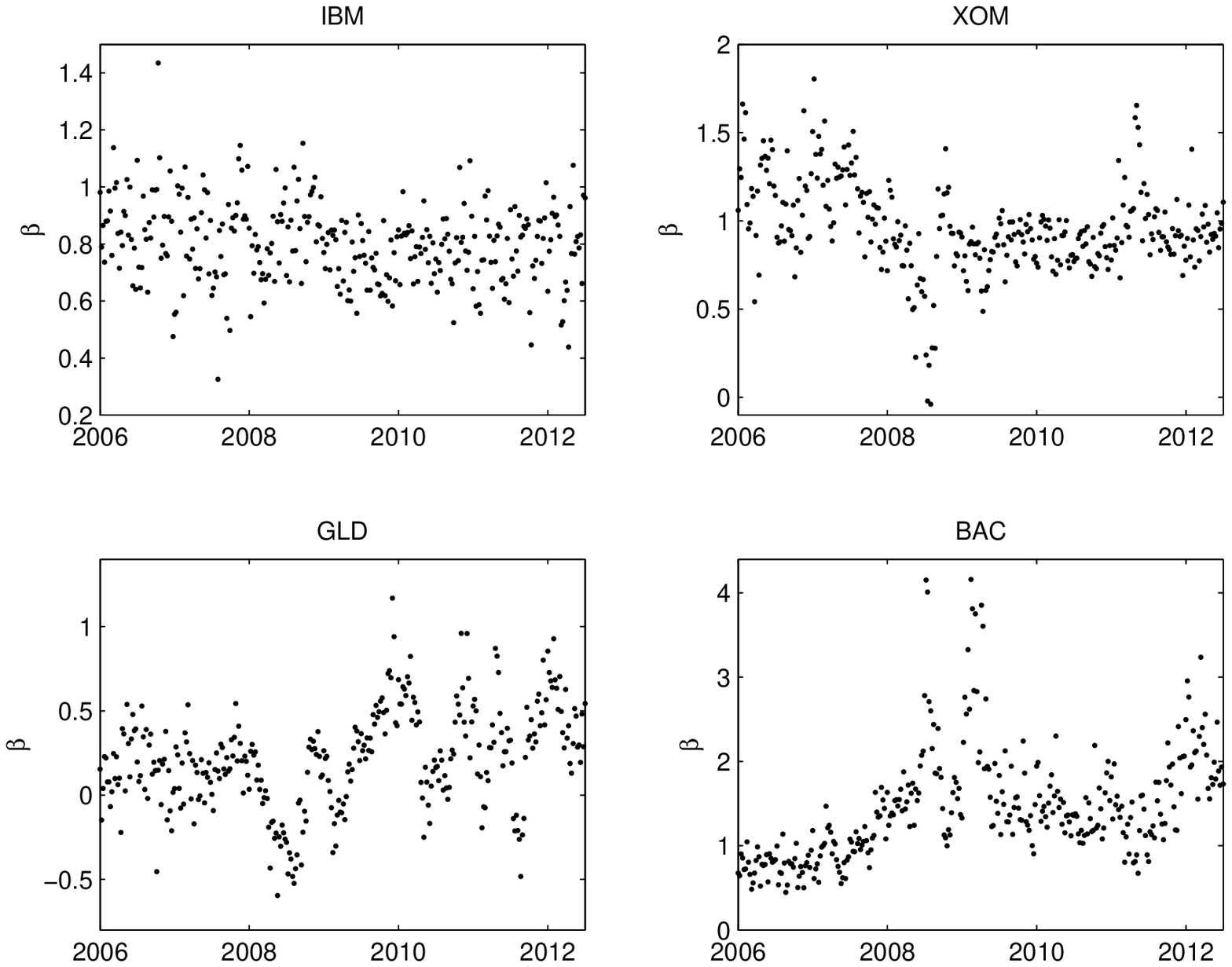}
\end{tabular}
\end{center}
\centering\caption{Estimated Betas Using Weekly Windows\label{fig:weekly_beta}. }
\end{figure}

Figure~\ref{fig:weekly_beta} shows time series of the weekly estimated market betas based on the pooled estimator in \refeq{beta_hat}. The contrasts are especially interesting when viewed in the context of the test results in Table~\ref{fig:weekly_beta} and also keeping in mind that the sample contains the most turbulent financial episode in many decades. Starting with IBM, we see from Figure~\ref{fig:weekly_beta} no significant pattern in the time series variation of its weekly market beta. Most of the weekly IBM market beta variation can be attributed to sampling error in its estimation. This is consistent with our test results in Table~\ref{tb:weekly_beta} and in particular the relatively low rejection rates for constancy of beta even over a time window of a month. Again consistent with our test results in Table~\ref{tb:weekly_beta} we see more time variation in the XOM market betas. Particularly noticeable is the period of July-August 2008 during which the market beta of XOM is quite low. Another interesting episode is that of May 2011 during which XOM's market beta was much higher than its average level.

Turning to GLD, we see a lot of variation in its sensitivity towards the market over the analyzed period. Perhaps not surprisingly, during bad times, such as the 2008 global financial crisis and the subsequent European monetary crises, GLD market beta is negative as during these periods gold serves the purpose of a hedging financial instrument. On the other hand, during normal times gold has little (positive) market sensitivity and it acts more like a pure commodity.  Finally for BAC, and completely in accordance with our results in Table~\ref{tb:weekly_beta}, we see very clear and persistent time variation. Over the period of 2006 until 2012, BAC market beta gradually increases from around 1 to around 2. Along the way of this gradual increase, we notice spikes in market beta around periods of crises such as the one in the Fall of 2008.

\section{Conclusion}\label{sec:concl}

We propose a nonparametric test for constant beta over a fixed interval of time from discrete observations of an asset and a risk factor with asymptotically vanishing distance between observations. The test is based on forming test statistics for constant beta over blocks with asymptotically increasing observations within them and shrinking time span and then summing them and scaling appropriately the resulting sum. The test is of self-normalized type which makes its limiting distribution under the null pivotal and independent from nuisance ``parameters'' such as the stochastic volatilities of the underlying processes.  We show asymptotic optimality for local nonparametric alternatives that are $\alpha$--H\"older regular. We find satisfactorily performance on simulated data in a Monte Carlo. In an empirical application we study the time window over which market betas of four different assets can be assumed to remain constant. 
\section{Proofs}\label{sec:proofs}

Throughout the proofs we will denote with $K$ a constant that does not depend on $n$ and the indices $i$ and $j$, but only on the characteristics of the multivariate process $(X,Y)$ and the powers involved in the estimates below, and further $K$ can change from line to line. We will further use the shorthand notation $\mathbb{E}_i^n\left(\cdot\right) = \mathbb{E}\left(\cdot|\mathcal{F}_{\frac{i-1}{n}}\right)$ and $\mathbb{P}_i^n\left(\cdot\right) = \mathbb{P}\left(\cdot|\mathcal{F}_{\frac{i-1}{n}}\right)$.

We start with some auxiliary notation to be used throughout the proofs. For arbitrary $b\in\mathbb{R}$, we denote
\begin{equation}\label{proof_1}
C_j^n(b) = \frac{n}{\sqrt{k_n}}\sum_{i=(j-1)k_n+1}^{jk_n}\Delta_i^nX^c(\Delta_i^nY^c-b\Delta_i^nX^c),~~V_j^n(b) = V_j^{(n,1)}V_j^{(n,2)}(b),\nonumber
\end{equation}
\begin{equation}\label{proof_2}
V_j^{(n,1)} = \frac{n}{k_n}\sum_{i =(j-1)k_n+1}^{jk_n}\left(\Delta_i^nX^c\right)^2,~~V_j^{(n,2)}(b) = \frac{n}{k_n}\sum_{i =(j-1)k_n+1}^{jk_n}\left(\Delta_i^nY^c-b\Delta_i^nX^c\right)^2,\nonumber
\end{equation}
\begin{equation}\label{proof_3}
T^n(b)= \frac{1}{\sqrt{2}}\sqrt{\frac{k_n}{n}}\sum_{j=2}^{\lfloor\frac{n}{k_n}\rfloor}T_j^n(b),~~T_j^n(b) =  \frac{\left(C_j^n(b)\right)^2-V_{j}^n(b)}{V_{j-1}^n(b)},\nonumber
\end{equation}
where recall $X^c$ and $Y^c$ are the continuous parts of the processes $X$ and $Y$. We also use the following shorthand notation
\begin{equation}\label{proof_4}
\widetilde{Y}_t^c = \int_0^t\widetilde{\alpha}_sds+\int_0^t\widetilde{\sigma}_sd\widetilde{W}_s,~~\widetilde{\alpha}_s = \alpha_s^Y-\beta_0\alpha_s^X,\nonumber
\end{equation}
and we further set
\begin{equation}\label{proof_5}
X_s^n = X_s^c-X_{\frac{i-1}{n}}^c,~~~\widetilde{Y}_s^n = \widetilde{Y}_s^c-\widetilde{Y}_{\frac{i-1}{n}}^c,~~\textrm{for $s\in\left[\frac{i-1}{n},\frac{i}{n}\right]$}.
\nonumber
\end{equation}
Finally, we denote
\begin{equation}\label{proof_6}
\mathcal{B}^n=\left\{|\widehat{\beta}_n-\overline{\beta}|\leq\delta n^{\iota-1/2}\right\},\nonumber
\end{equation}
for some arbitrary small positive numbers $\iota>0$ and $\delta>0$.

\subsection{Localization}\label{subsec:proofs_loc}
We will proof the results under the following stronger assumption:

\noindent \textbf{Assumption SA.} For the process defined in \refeq{eq:XY}-\refeq{eq:vol} we have:
\begin{itemize}
\item [(a)] \textit{$|\sigma_t|$ and $|\widetilde{\sigma}_{t}|$ are uniformly bounded from below and above; }
\item [(b)] \textit{ $\beta$, $\alpha^{\sigma}$, $\alpha^{\widetilde{\sigma}}$, $\gamma^{\sigma}$, $\gamma^{\widetilde{\sigma}}$, $\widetilde{\gamma}^{\sigma}$, $\widetilde{\gamma}^{\widetilde{\sigma}}$, $\gamma'$ and $\gamma^{''}$ are bounded; $\delta^X$, $\delta^Y$, $\delta^{\sigma}$ and $\delta^{\widetilde{\sigma}}$ are bounded;}
\item [(c)] \textit{the coefficients in the \Ito semimartingale representations of $\alpha^X$ and $\alpha^Y$ are bounded;}
\item [(d)] \textit{$|\delta^{X}(t,x)|+|\delta^{Y}(t,x)|\leq \gamma^{(1)}(x)$ for all $t\leq 1$ with $\int_{\mathbb{R}}|\gamma^{(1)}(x)|^{r}dx<\infty$ for some $r\in(0,1)$; $|\delta^{\sigma}(t,x)|+|\delta^{\widetilde{\sigma}}(t,x)|\leq \gamma^{(2)}(x)$ for all $t\leq 1$ with $\int_{\mathbb{R}}|\gamma^{(2)}(x)|^{2}dx<\infty$;}
\end{itemize}
Extending the results to the case when only the weaker assumption A holds follows from standard localization procedure as in Lemma 4.4.9 of \cite{JP12}.

\subsection{Preliminary results}
In this section we will always assume (without further mention) that $k_n\rightarrow\infty$ and $k_n/n\rightarrow 0$.
\begin{lema}\label{lema_mean}
Under assumption SA and $\beta_t=\beta_0$ for $t\in[0,1]$, we have
\begin{equation}\label{lema_mean_1}
\left| \mathbb{E}_{\frac{(j-1)k_n}{n}}\left( C_j^n(\beta_0)^2 - V_j^n(\beta_0) \right) \right|\leq K\left(\frac{k_n}{n}\bigvee\frac{1}{\sqrt{n}}\right).
\end{equation}
\end{lema}
\textbf{Proof of Lemma~\ref{lema_mean}.}
First, we derive some bounds for $C_j^n(\beta_0)^2$. Note that for $\beta_t=\beta_0$ on $t\in[0,1]$, we have $C_j^n(\beta_0) = \frac{n}{\sqrt{k_n}}\sum_{i=(j-1)k_n+1}^{jk_n}\Delta_i^nX^c\Delta_i^n\widetilde{Y}^c$. Applying \Ito formula and since $W$ and $\widetilde{W}$ are orthogonal, we have
\begin{equation}\label{lema_mean_proof_1}
\mathbb{E}_{i-1}^n(\Delta_i^nX^c\Delta_i^n\widetilde{Y}^c) = \mathbb{E}_{i-1}^n\left(\int_{\frac{i-1}{n}}^{\frac{i}{n}}X_s^n\widetilde{\alpha}_sds+\int_{\frac{i-1}{n}}^{\frac{i}{n}}\widetilde{Y}_s^n\alpha_s^Xds\right).
\end{equation}
Next, given the \Ito semimartingale assumption for the processes $\alpha^X$ and $\alpha^Y$, as well as an application of Cauchy-Schwarz and Burkholder-Davis-Gundy inequalities, we have
\begin{equation}\label{lema_mean_proof_2}
\left|\mathbb{E}_{i-1}^n\left(\int_{\frac{i-1}{n}}^{\frac{i}{n}}X_s^n(\widetilde{\alpha}_s-\widetilde{\alpha}_{\frac{i-1}{n}})ds+\int_{\frac{i-1}{n}}^{\frac{i}{n}}
\widetilde{Y}_s^n(\alpha_s^X-\alpha_{\frac{i-1}{n}}^X)ds\right)\right|\leq \frac{K}{n^2}.\nonumber
\end{equation}
From here, using the definition of the processes $X^n$ and $\widetilde{Y}^n$, we have altogether
\begin{equation}\label{lema_mean_proof_3}
|\mathbb{E}_{i-1}^n(\Delta_i^nX^c\Delta_i^n\widetilde{Y}^c)|\leq \frac{K}{n^2}.
\end{equation}
Next, using \Ito formula we have
\begin{equation}\label{lema_mean_proof_4}
(\Delta_i^nX^c)^2(\Delta_i^n\widetilde{Y}^c)^2 = \left(2\int_{\frac{i-1}{n}}^{\frac{i}{n}}X_s^ndX_s^c+ \int_{\frac{i-1}{n}}^{\frac{i}{n}}\sigma_s^2ds \right)
\left(2\int_{\frac{i-1}{n}}^{\frac{i}{n}}\widetilde{Y}_s^nd\widetilde{Y}_s^c+ \int_{\frac{i-1}{n}}^{\frac{i}{n}}\widetilde{\sigma}_s^2ds \right).\nonumber
\end{equation}
Applying \Ito formula, Cauchy-Schwarz and Burkholder-Davis-Gundy inequalities, and using the independence of $W$ and $\widetilde{W}$, we get
\begin{equation}\label{lema_mean_proof_5}
\left|\mathbb{E}_{i-1}^n\left(\int_{\frac{i-1}{n}}^{\frac{i}{n}}X_s^ndX_s^c\int_{\frac{i-1}{n}}^{\frac{i}{n}}\widetilde{Y}_s^nd\widetilde{Y}_s^c\right) \right|\leq \frac{K}{n^2\sqrt{n}}.
\end{equation}
Further using the \Ito semimartingale assumption for $\sigma$, Cauchy-Schwarz and Burkholder-Davis-Gundy inequalities, we have
\begin{equation}\label{lema_mean_proof_6}
\left|\mathbb{E}_{i-1}^n\left(\int_{\frac{i-1}{n}}^{\frac{i}{n}}X_s^ndX_s^c\int_{\frac{i-1}{n}}^{\frac{i}{n}}(\widetilde{\sigma}_s^2 -\widetilde{\sigma}_{\frac{i-1}{n}}^2 )ds   \right)\right|+
\left|\mathbb{E}_{i-1}^n\left(\int_{\frac{i-1}{n}}^{\frac{i}{n}}\widetilde{Y}_s^nd\widetilde{Y}_s^c\int_{\frac{i-1}{n}}^{\frac{i}{n}}(\sigma_s^2 -\sigma_{\frac{i-1}{n}}^2 )ds   \right)\right|\leq \frac{K}{n^2\sqrt{n}}.
\end{equation}
Finally using the definition of $X^n$ and $\widetilde{Y}^n$, exactly as in \refeq{lema_mean_proof_1} and \refeq{lema_mean_proof_3} above, we get
\begin{equation}\label{lema_mean_proof_7}
\left|\mathbb{E}_{i-1}^n\left(\int_{\frac{i-1}{n}}^{\frac{i}{n}}X_s^ndX_s^c  \right)\right|+
\left|\mathbb{E}_{i-1}^n\left(\int_{\frac{i-1}{n}}^{\frac{i}{n}}\widetilde{Y}_s^nd\widetilde{Y}_s^c   \right)\right|\leq \frac{K}{n^2}.
\end{equation}
Using the bounds in \refeq{lema_mean_proof_3}-\refeq{lema_mean_proof_7}, we get
\begin{equation}\label{lema_mean_proof_8}
\mathbb{E}_{\frac{(j-1)k_n}{n}}\left(C_j^n(\beta^n)\right)^2 = \frac{n^2}{k_n}\sum_{i=(j-1)k_n+1}^{jk_n}\mathbb{E}_{\frac{(j-1)k_n}{n}}
\left(\int_{\frac{i-1}{n}}^{\frac in}\sigma_s^2ds\int_{\frac{i-1}{n}}^{\frac in}\widetilde{\sigma}_s^2ds\right)+\widetilde{R}_j^{(n,1)},\nonumber
\end{equation}
where
\begin{equation}\label{lema_mean_proof_9}
\left|\widetilde{R}_j^{(n,1)}\right|\leq K\left(\frac{k_n}{n}\bigvee\frac{1}{\sqrt{n}}\right).
\end{equation}
We turn next to $V_j^n(\beta_0)$. Using \Ito formula, we can write
\begin{equation}\label{lema_mean_proof_10}
\left\{\begin{array}{l}
V_j^{(n,1)} = \frac{n}{k_n}\left(2\int_{\frac{(j-1)k_n}{n}}^{\frac{jk_n}{n}}X_s^ndX_s^c+ \int_{\frac{(j-1)k_n}{n}}^{\frac{jk_n}{n}}\sigma_s^2ds \right),\\
V_j^{(n,2)}(\beta_0) = \frac{n}{k_n}\left(2\int_{\frac{(j-1)k_n}{n}}^{\frac{jk_n}{n}}\widetilde{Y}_s^nd\widetilde{Y}_s^c+ \int_{\frac{(j-1)k_n}{n}}^{\frac{jk_n}{n}}\widetilde{\sigma}_s^2ds \right).
\end{array}\right.
\end{equation}
From here, using similar bounds to the ones derived in \refeq{lema_mean_proof_5}-\refeq{lema_mean_proof_7}, we get
\begin{equation}\label{lema_mean_proof_11}
\mathbb{E}_{\frac{(j-1)k_n}{n}}\left(V_j^n(\beta_0)\right) = \frac{n^2}{k_n^2}\mathbb{E}_{\frac{(j-1)k_n}{n}}
\left(\int_{\frac{(j-1)k_n}{n}}^{\frac{jk_n}{n}}\sigma_s^2ds\int_{\frac{(j-1)k_n}{n}}^{\frac{jk_n}{n}}\widetilde{\sigma}_s^2ds\right)
+\widetilde{R}_j^{(n,2)},\nonumber
\end{equation}
where
\begin{equation}\label{lema_mean_proof_12}
\left|\widetilde{R}_j^{(n,2)}\right|\leq \frac{K}{\sqrt{n}}.
\end{equation}
Given the bounds for the conditional expectations of the residual terms $\widetilde{R}_j^{(n,1)}$ and $\widetilde{R}_j^{(n,2)}$, we are left with the difference
\[\frac{n^2}{k_n}\sum_{i=(j-1)k_n+1}^{jk_n}\mathbb{E}_{\frac{(j-1)k_n}{n}}
\left(\int_{\frac{i-1}{n}}^{\frac{i}{n}}\sigma_s^2ds\int_{\frac{i-1}{n}}^{\frac{i}{n}}\widetilde{\sigma}_s^2ds\right) - \frac{n^2}{k_n^2}\mathbb{E}_{\frac{(j-1)k_n}{n}}
\left(\int_{\frac{(j-1)k_n}{n}}^{\frac{jk_n}{n}}\sigma_s^2ds\int_{\frac{(j-1)k_n}{n}}^{\frac{jk_n}{n}}\widetilde{\sigma}_s^2ds\right).\]
Using the \Ito semimartingale representation of $\sigma^2$ and $\widetilde{\sigma}^2$ in \refeq{eq:vol} and Cauchy-Schwarz inequality, we have
\begin{equation}\label{lema_mean_proof_13}
|\mathbb{E}\left(\sigma_t^2-\sigma_u^2|\mathcal{F}_u\right)|+|\mathbb{E}\left(\widetilde{\sigma}_t^2-\widetilde{\sigma}_u^2|\mathcal{F}_u\right)|\leq K|t-u|,~~u\leq t,
\nonumber
\end{equation}
\begin{equation}\label{lema_mean_proof_14}
|\mathbb{E}\left((\sigma_t^2-\sigma_u^2)(\widetilde{\sigma}_t^2-\widetilde{\sigma}_u^2)|\mathcal{F}_u\right)|\leq K|t-u|,~~u\leq t.\nonumber
\end{equation}
Using these inequalities as well as the algebraic identity
\begin{equation}\label{eq:algidentity}
x_1y_1-x_2y_2 = (x_1-x_2)(y_1-y_2)+x_2(y_1-y_2)+(x_1-x_2)y_2, \text{ for any real }x_1, x_2, y_1,y_2,
\end{equation}
 we have
\begin{equation}\label{lema_mean_proof_15}
\left|\mathbb{E}_{\frac{(j-1)k_n}{n}}\left( \int_{\frac{i-1}{n}}^{\frac{i}{n}}\sigma_s^2ds\int_{\frac{i-1}{n}}^{\frac{i}{n}}\widetilde{\sigma}_s^2ds - \frac{1}{n^2}\sigma_{\frac{(j-1)k_n}{n}}^2\widetilde{\sigma}_{\frac{(j-1)k_n}{n}}^2  \right)\right|\leq K\frac{k_n}{n^3},~~i=(j-1)k_n+1,..,jk_n,\nonumber
\end{equation}
\begin{equation}\label{lema_mean_proof_16}
\left|\mathbb{E}_{\frac{(j-1)k_n}{n}}\left( \int_{\frac{(j-1)k_n}{n}}^{\frac{jk_n}{n}}\sigma_s^2ds\int_{\frac{(j-1)k_n}{n}}^{\frac{jk_n}{n}}\widetilde{\sigma}_s^2ds - \frac{k_n^2}{n^2}\sigma_{\frac{(j-1)k_n}{n}}^2\widetilde{\sigma}_{\frac{(j-1)k_n}{n}}^2  \right)\right|\leq K\left(\frac{k_n}{n}\right)^3.\nonumber
\end{equation}
Combining these results with the bounds in \refeq{lema_mean_proof_9} and \refeq{lema_mean_proof_12}, we get the result to be proved.
\qedd

\begin{lema}\label{lema_mean_2}
Under assumption SA and $\beta_t=\beta_0$ for $t\in[0,1]$, we have
\begin{equation}\label{lema_mean_2_1}
\mathbb{E}_{\frac{(j-1)k_n}{n}}\left|V_j^n(\beta_0) -\frac{n^2}{k_n^2}\int_{\frac{(j-1)k_n}{n}}^{\frac{jk_n}{n}}\sigma_{s}^2ds\int_{\frac{(j-1)k_n}{n}}^{\frac{jk_n}{n}}\widetilde{\sigma}_{s}^2ds\right|^2\leq \frac{K}{k_n}.
\end{equation}
\end{lema}

\textbf{Proof of Lemma~\ref{lema_mean_2}.}
We make use of \eqref{eq:algidentity} as well as the boundedness of the processes $\sigma^2$ and $\widetilde{\sigma}^2$, to bound
\begin{equation}\label{lema_mean_2_proof_1}
\begin{split}
&\mathbb{E}_{\frac{(j-1)k_n}{n}}\left|V_j^n(\beta_0) -\frac{n^2}{k_n^2}\int_{\frac{(j-1)k_n}{n}}^{\frac{jk_n}{n}}\sigma_{s}^2ds\int_{\frac{(j-1)k_n}{n}}^{\frac{jk_n}{n}}\widetilde{\sigma}_{s}^2ds\right|^2\\&~~\leq K\mathbb{E}_{\frac{(j-1)k_n}{n}}\left|V_j^{(n,1)} -\frac{n}{k_n}\int_{\frac{(j-1)k_n}{n}}^{\frac{jk_n}{n}}\sigma_{s}^2ds\right|^2+K\mathbb{E}_{\frac{(j-1)k_n}{n}}\left|V_j^{(n,2)}(\beta_0) -\frac{n}{k_n}\int_{\frac{(j-1)k_n}{n}}^{\frac{jk_n}{n}}\widetilde{\sigma}_{s}^2ds\right|^2\\&~~~~+K\mathbb{E}_{\frac{(j-1)k_n}{n}}\left|V_j^{(n,1)} -\frac{n}{k_n}\int_{\frac{(j-1)k_n}{n}}^{\frac{jk_n}{n}}\sigma_{s}^2ds\right|^4+K\mathbb{E}_{\frac{(j-1)k_n}{n}}\left|V_j^{(n,2)}(\beta_0) -\frac{n}{k_n}\int_{\frac{(j-1)k_n}{n}}^{\frac{jk_n}{n}}\widetilde{\sigma}_{s}^2ds\right|^4.\nonumber
\end{split}
\end{equation}
Using the decomposition of $V_j^{(n,1)}$ and $V_j^{(n,2)}(\beta_0)$ in \refeq{lema_mean_proof_10} and applying the Burkholder-Davis-Gundy inequality we get the result to be proved. \qedd

\begin{lema}\label{lema_var}
Under assumption SA and $\beta_t=\beta_0$ for $t\in[0,1]$, we have
\begin{equation}\label{lema_var_1}
\left|\mathbb{E}_{\frac{(j-1)k_n}{n}}\left( C_j^n(\beta_0)^2 - V_j^n(\beta_0) \right)^2- 2\sigma_{\frac{(j-1)k_n}{n}}^4\widetilde{\sigma}_{\frac{(j-1)k_n}{n}}^4\right|\leq K\left[\left(\frac{k_n}{n}\right)^{1/2-\iota}\bigvee\frac{1}{\sqrt{k_n}}\right],~~~\forall \iota>0.
\end{equation}
\end{lema}
\textbf{Proof of Lemma~\ref{lema_var}.}
We first denote the analogues of $C_j^n(\beta_0)$ and $V_j^n(\beta_0)$, with $\sigma_s$ and $\widetilde{\sigma}_s$ kept at their values at the beginning of the block, as
\begin{equation}\label{lema_var_proof_1}
\begin{split}
\overline{C}_j^n &= \frac{n}{\sqrt{k_n}}\sigma_{\frac{(j-1)k_n}{n}}\widetilde{\sigma}_{\frac{(j-1)k_n}{n}}\sum_{i=(j-1)k_n+1}^{jk_n}\Delta_i^nW\Delta_i^n\widetilde{W},~\overline{V}_j^n = \overline{V}_j^{(n,1)}\overline{V}_j^{(n,2)},\\
\overline{V}_j^{(n,1)} &= \frac{n}{k_n}\sigma_{\frac{(j-1)k_n}{n}}^2\sum_{i=(j-1)k_n+1}^{jk_n}(\Delta_i^nW)^2,~~\overline{V}_j^{(n,2)} =\frac{n}{k_n}
\widetilde{\sigma}_{\frac{(j-1)k_n}{n}}^2\sum_{i=(j-1)k_n+1}^{jk_n}(\Delta_i^n\widetilde{W})^2.\nonumber
\end{split}
\end{equation}
Using the independence of $W$ and $\widetilde{W}$ and Burkholder-Davis-Gundy inequality for discrete martingales, we have
\begin{equation}\label{lema_var_proof_3}
\mathbb{E}_{\frac{(j-1)k_n}{n}}\left|\overline{C}_j^n\right|^p\leq K,~~~\mathbb{E}_{\frac{(j-1)k_n}{n}}\left|\overline{V}_j^{(n,1)}\right|^p\leq K,~~~\mathbb{E}_{\frac{(j-1)k_n}{n}}\left|\overline{V}_j^{(n,2)}\right|^p\leq K,~~~\forall p\geq 2.
\end{equation}
Using the algebraic identity $x^2-y^2 = (x-y)^2+2y(x-y)$ as well as Cauchy-Schwarz inequality and \refeq{lema_var_proof_3}, we can write
\begin{equation}\label{lema_var_proof_2}
\begin{split}
&\left|\mathbb{E}_{\frac{(j-1)k_n}{n}}\left( C_j^n(\beta_0)^2 - V_j^n(\beta_0) \right)^2 - \mathbb{E}_{\frac{(j-1)k_n}{n}}\left( (\overline{C}_j^n)^2 - \overline{V}_j^n \right)^2\right|\\&~~~~~~~~~~~~~~~~~~~~~~~~\leq
K \chi\left(\mathbb{E}_{\frac{(j-1)k_n}{n}}\left( C_j^n(\beta_0)^2 - (\overline{C}_j^n)^2 \right)^2\right) +
K \chi\left(\mathbb{E}_{\frac{(j-1)k_n}{n}}\left( V_j^n(\beta_0) - \overline{V}_j^n \right)^2\right),
\end{split}
\end{equation}
for $\chi(u) = u\vee\sqrt{u}$. Next, applying the \Ito formula and using the \Ito semimartingale assumption for the process $\alpha^X$, we have
\begin{equation}\label{lema_var_proof_4}
\left|\mathbb{E}_{i-1}^n\left(\Delta_i^nX^c-\sigma_{\frac{(j-1)k_n}{n}}\Delta_i^nW\right)\Delta_i^n\widetilde{Y}^c \right|\leq K\frac{\sqrt{k_n}}{n^2}.\nonumber
\end{equation}
From here, using Burkholder-Davis-Gundy inequality for discrete martingales, together with our assumption for $\sigma$ being \Ito semimartingale, \Ito formula and the independence of $W$ from $\widetilde{W}$, as well as \Holder~ inequality, we have for every $p\geq 1$ and any $\iota>0$
\begin{equation}\label{lema_var_proof_5}
\begin{split}
&\mathbb{E}_{\frac{(j-1)k_n}{n}}\left|\frac{n}{\sqrt{k_n}}\sum_{i=(j-1)k_n+1}^{jk_n}\left(\Delta_i^nX^c-\sigma_{\frac{(j-1)k_n}{n}}\Delta_i^nW\right)
\Delta_i^n\widetilde{Y}^c \right|^p\leq K\left(\frac{k_n}{n}\right)^{(p/2)\wedge 1-\iota}.\nonumber
\end{split}
\end{equation}
Similar analysis implies
\begin{equation}\label{lema_var_proof_6}
\begin{split}
&\mathbb{E}_{\frac{(j-1)k_n}{n}}\left|\frac{n}{\sqrt{k_n}}\sum_{i=(j-1)k_n+1}^{jk_n}\Delta_i^nW\left(\Delta_i^n\widetilde{Y}^c-\widetilde{\sigma}_{\frac{(j-1)k_n}{n}}
\Delta_i^n\widetilde{W}\right) \right|^p\leq K\left(\frac{k_n}{n}\right)^{(p/2)\wedge 1-\iota},\nonumber
\end{split}
\end{equation}
and therefore
\begin{equation}\label{lema_var_proof_7}
\mathbb{E}_{\frac{(j-1)k_n}{n}}\left|C_j^n(\beta_0) - \overline{C}_j^n \right|^p\leq K\left(\frac{k_n}{n}\right)^{(p/2)\wedge 1-\iota}.
\end{equation}
Combining this result with the bound in \refeq{lema_var_proof_3}, together with Cauchy-Schwarz inequality, we get
\begin{equation}\label{lema_var_proof_8}
\mathbb{E}_{\frac{(j-1)k_n}{n}}\left( C_j^n(\beta_0)^2 - (\overline{C}_j^n)^2 \right)^2\leq K\left(\frac{k_n}{n}\right)^{1-\iota}.
\end{equation}

We next bound $\mathbb{E}_{\frac{(j-1)k_n}{n}}\left( V_j^n(\beta_0) - \overline{V}_j^n \right)^2$. Using inequality in means, i.e., $\left|\frac{\sum_{i=1}^nx_i}{n}\right|^p\leq \frac{\sum_{i=1}^n|x_i|^p}{n}$ for any $p\geq 1$ and any real $\{x_i\}_{i=1,..,n}$, we first have for $p\geq 1$
\begin{small}
\begin{equation}\label{lema_var_proof_9}
\mathbb{E}_{\frac{(j-1)k_n}{n}}\left(\frac{n}{k_n}\sum_{i=(j-1)k_n+1}^{jk_n}\left[(\Delta_i^nX^c)^2+(\Delta_i^n\widetilde{Y}^c)^2
+\sigma_{\frac{(j-1)k_n}{n}}^2(\Delta_i^nW)^2 + \widetilde{\sigma}_{\frac{(j-1)k_n}{n}}^2(\Delta_i^n\widetilde{W})^2\right]  \right)^p\leq K.\nonumber
\end{equation}
\end{small}
Using the algebraic identity $x^2-y^2 = (x-y)^2+2y(x-y)$, the Burkholder-Davis-Gundy inequality for discrete martingales, \Ito formula for the function $f(x,y) = x y$, our assumption for $\sigma$ being \Ito semimartingale, we have for any $\iota>0$
\begin{equation}\label{lema_var_proof_10}
\mathbb{E}_{\frac{(j-1)k_n}{n}}\left|\frac{n}{k_n}\sum_{i=(j-1)k_n+1}^{jk_n}\left[(\Delta_i^nX^c)^2-\sigma_{\frac{(j-1)k_n}{n}}^2(\Delta_i^nW)^2\right]\right|^p\leq K\left(\frac{k_n}{n}\right)^{1-\iota},~~~\forall p\geq 2,\nonumber
\end{equation}
and similarly
\begin{equation}\label{lema_var_proof_11}
\mathbb{E}_{\frac{(j-1)k_n}{n}}\left|\frac{n}{k_n}\sum_{i=(j-1)k_n+1}^{jk_n}\left[(\Delta_i^n\widetilde{Y}^c)^2
-\widetilde{\sigma}_{\frac{(j-1)k_n}{n}}^2(\Delta_i^n\widetilde{W})^2\right]\right|^p\leq K\left(\frac{k_n}{n}\right)^{1-\iota},~~~\forall p\geq 2.\nonumber
\end{equation}
Using the above inequalities and \Holder~ inequality, we have
\begin{equation}\label{lema_var_proof_12}
\mathbb{E}_{\frac{(j-1)k_n}{n}}\left( V_j^n(\beta_0) - \overline{V}_j^n \right)^2\leq K\left(\frac{k_n}{n}\right)^{1-\iota},~~~~\forall \iota>0.
\end{equation}
Altogether, combining the bounds in \refeq{lema_var_proof_2}, \refeq{lema_var_proof_8} and \refeq{lema_var_proof_12}, we get
\begin{equation}\label{lema_var_proof_13}
\begin{split}
\left|\mathbb{E}_{\frac{(j-1)k_n}{n}}\left( C_j^n(\beta_0)^2 - V_j^n(\beta_0) \right)^2 - \mathbb{E}_{\frac{(j-1)k_n}{n}}\left( (\overline{C}_j^n)^2 - \overline{V}_j^n \right)^2\right|\leq
K \left(\frac{k_n}{n}\right)^{1/2-\iota},~~~~\forall \iota>0.
\end{split}
\end{equation}
We are thus left with $\mathbb{E}_{\frac{(j-1)k_n}{n}}\left((\overline{C}_j^n)^2-\overline{V}_j^n\right)^2$. First, using finite sample distribution results for regressions with normally distributed errors, see e.g., \cite{Hayashi}, we have
\begin{equation}\label{lema_var_proof_14}
\frac{\frac{(\overline{C}_j^n)^2}{\overline{V}_j^n}}{1-\frac{1}{k_n}\frac{(\overline{C}_j^n)^2}{\overline{V}_j^n}} ~\stackrel{d}{=}~\frac{k_n}{k_n-1}t_{k_n-1}^2~~~\Longrightarrow
\frac{(\overline{C}_j^n)^2}{\overline{V}_j^n}~\stackrel{d}{=}~\frac{t_{k_n-1}^2}{1+\frac{1}{k_n}(t_{k_n-1}^2-1)},\nonumber
\end{equation}
where $t_k$ denotes a random variable, having a $t$-distribution with $k$ degrees of freedom. Therefore, for $k_n>9$ (so that the $t_{k_n-1}$-distribution has finite eight moment), using the moments of the $t$-distribution, we have
\begin{equation}\label{lema_var_proof_15}
\left|\mathbb{E}_{\frac{(j-1)k_n}{n}}\left(\frac{(\overline{C}_j^n)^2}{\overline{V}_j^n}-1\right)^2-2\right|\leq \frac{K}{k_n}.\nonumber
\end{equation}
Second, using Burkholder-Davis-Gundy inequality for discrete martingales, we have
\begin{equation}\label{lema_var_proof_16}
\mathbb{E}_{\frac{(j-1)k_n}{n}}\left|\left( \overline{V}_j^n \right)^2 - \sigma_{\frac{(j-1)k_n}{n}}^4\widetilde{\sigma}_{\frac{(j-1)k_n}{n}}^4\right|^p\leq \frac{K}{k_n^{p/2}},~~~\forall p\geq 2.\nonumber
\end{equation}
Combining the above bounds, and using \Holder~ inequality, we have for $k_n>9$
\begin{equation}\label{lema_var_proof_17}
\left| \mathbb{E}_{\frac{(j-1)k_n}{n}}\left((\overline{C}_j^n)^2-\overline{V}_j^n\right)^2 - 2 \sigma_{\frac{(j-1)k_n}{n}}^4\widetilde{\sigma}_{\frac{(j-1)k_n}{n}}^4 \right|\leq \frac{K}{\sqrt{k_n}}.
\end{equation}
The result of the lemma then follows from \refeq{lema_var_proof_13} and \refeq{lema_var_proof_17}.
\qedd
\begin{lema}\label{lema_powers}
Under assumption SA and $\beta_t=\beta_0$ for $t\in[0,1]$, for any constant $\alpha>0$, we have
\begin{equation}\label{lema_powers_1}
\mathbb{E}_{\frac{(j-1)k_n}{n}}\left|C_j^n(\beta_0)\right|^p+\mathbb{E}_{\frac{(j-1)k_n}{n}}\left|V_j^{(n,1)}\right|^p
+\mathbb{E}_{\frac{(j-1)k_n}{n}}\left|V_j^{(n,2)}(\beta_0)\right|^p\leq K,~~\forall p\geq 1,
\end{equation}
\begin{equation}\label{lema_powers_2}
\mathbb{P}\left(\left| V_{j}^{(n,1)} - \frac{n}{k_n}\int_{\frac{(j-1)k_n}{n}}^{\frac{jk_n}{n}}\sigma_{s}^2ds \right|\geq \alpha \frac{n}{k_n} \int_{\frac{(j-1)k_n}{n}}^{\frac{jk_n}{n}}\sigma_{s}^2ds\right)\leq \frac{K}{k_n^{p/2}},~~\forall p\geq 1,
\end{equation}
\begin{equation}\label{lema_powers_3}
\mathbb{P}\left(\left| V_{j}^{(n,2)}(\beta_0) - \frac{n}{k_n}\int_{\frac{(j-1)k_n}{n}}^{\frac{jk_n}{n}}\widetilde{\sigma}_{s}^2ds \right|\geq \alpha \frac{n}{k_n} \int_{\frac{(j-1)k_n}{n}}^{\frac{jk_n}{n}}\widetilde{\sigma}_{s}^2ds\right)\leq \frac{K}{k_n^{p/2}},~~\forall p\geq 1,
\end{equation}
where the constant $K$ in the above bounds depends on the constant $\alpha$.
\end{lema}
\textbf{Proof of Lemma~\ref{lema_powers}.}
Using Burkholder-Davis-Gundy inequality for discrete martingales (and \Holder~inequality when $p<2$), we have
\begin{equation}\label{lema_powers_proof_1}
\mathbb{E}_{\frac{(j-1)k_n}{n}}\left|
\frac{n}{\sqrt{k_n}}\sum_{i=(j-1)k_n+1}^{jk_n}\left[(\Delta_i^nX^c\Delta_i^n\widetilde{Y}^c)-\mathbb{E}_{i-1}^n(\Delta_i^nX^c\Delta_i^n\widetilde{Y}^c)\right]\right|^p\leq K,~~~p\geq 1.\nonumber
\end{equation}
Using this bound together with the bound for $\mathbb{E}_{i-1}^n(\Delta_i^nX^c\Delta_i^n\widetilde{Y}^c)$ in \refeq{lema_mean_proof_3} (note that $\sqrt{k_n}/n\rightarrow 0$), we get the bound for $C_j^n(\beta_0)$ in \refeq{lema_powers_1}. The bounds for $V_j^{(n,1)}$ and $V_j^{(n,2)}(\beta_0)$ in \refeq{lema_powers_1} follow from inequality in means.

Next using the decomposition of $V_j^{(n,1)}$ and $V_j^{(n,2)}(\beta_0)$ in \refeq{lema_mean_proof_10}, we have by an application of Burkholder-Davis-Gundy inequality
\begin{equation}\label{lema_powers_proof_2}
\left\{\begin{array}{l}
\mathbb{E}_{\frac{(j-1)k_n}{n}}\left|V_j^{(n,1)}-\frac{n}{k_n}\int_{\frac{(j-1)k_n}{n}}^{\frac{jk_n}{n}}\sigma_{s}^2ds\right|^p\leq \frac{K}{k_n^{p/2}},~~\forall p\geq 1,
\\ \mathbb{E}_{\frac{(j-1)k_n}{n}}\left|V_j^{(n,2)}(\beta_0)-\frac{n}{k_n}\int_{\frac{(j-1)k_n}{n}}^{\frac{jk_n}{n}}\widetilde{\sigma}_{s}^2ds\right|^p\leq \frac{K}{k_n^{p/2}},~~\forall p\geq 1.\end{array}\right.\nonumber
\end{equation}
From here, using the boundedness of the processes $|\sigma|$ and $|\widetilde{\sigma}|$, both from below and above, we get the bounds in \refeq{lema_powers_2} and \refeq{lema_powers_3}.
\qedd

\begin{lema}\label{lema:mom_altern}
Under Assumption SA for any constant $\alpha>0$ and provided $n^{2\iota-1}k_n\rightarrow 0$, for $\iota$ being the constant in the definition of the set $\mathcal{B}^n$, we have
\begin{small}
\begin{equation}\label{lema:mom_altern_1}
\mathbb{E}_{\frac{(j-1)k_n}{n}}\left\{\left(\frac{1}{\sqrt{k_n}}C_j^n(\widehat{\beta}_n)-
\frac{n}{k_n}\int_{\frac{(j-1)k_n}{n}}^{\frac{jk_n}{n}}(\beta_{s}-\overline{\beta})\sigma_{s}^2ds\right)^2
1_{\{\mathcal{B}^n\}}\right\}\leq \frac{K}{k_n},
\end{equation}
\begin{equation}\label{lema:mom_altern_2}
\mathbb{E}_{\frac{(j-1)k_n}{n}}\left\{\left(V_j^{(n,1)}-\frac{n}{k_n}\int_{\frac{(j-1)k_n}{n}}^{\frac{jk_n}{n}}\sigma_{s}^2ds\right)^2\right\}\leq \frac{K}{k_n},
\end{equation}
\begin{equation}\label{lema:mom_altern_3}
\mathbb{E}_{\frac{(j-1)k_n}{n}}\left\{\left(V_j^{(n,2)}(\widehat{\beta}_n)-\frac{n}{k_n}\int_{\frac{(j-1)k_n}{n}}^{\frac{jk_n}{n}}((\beta_{s}-\overline{\beta})^2\sigma_{s}^2
+\widetilde{\sigma}^2_{s})ds\right)^21_{\{\mathcal{B}^n\}}\right\}\leq \frac{K}{k_n},
\end{equation}
\begin{equation}\label{lema:mom_altern_4}
\mathbb{P}\left(\left|V_j^{(n,1)}-\frac{n}{k_n}\int_{\frac{(j-1)k_n}{n}}^{\frac{jk_n}{n}}\sigma_s^2ds\right|
1_{\{\mathcal{B}^n\}}>\alpha\frac{n}{k_n}\int_{\frac{(j-1)k_n}{n}}^{\frac{jk_n}{n}}\sigma_s^2ds\right)\leq \frac{K}{k_n^{p/2}},~~\forall p\geq1,
\end{equation}
\begin{equation}\label{lema:mom_altern_5}
\mathbb{P}\left(\left|V_j^{(n,2)}(\widehat{\beta}_n)-\frac{n}{k_n}\int_{\frac{(j-1)k_n}{n}}^{\frac{jk_n}{n}}((\beta_s-\widehat{\beta}_n)^2\sigma_s^2+\widetilde{\sigma}_s^2)
ds\right|1_{\{\mathcal{B}^n\}}>\alpha\frac{n}{k_n}\int_{\frac{(j-1)k_n}{n}}^{\frac{jk_n}{n}}\widetilde{\sigma}_s^2ds\right)\leq \frac{K}{k_n^{p/2}},~~\forall p\geq 1.
\end{equation}
\end{small}
\end{lema}
\noindent \textbf{Proof of Lemma~\ref{lema:mom_altern}.}
Using \Ito formula, we have
\begin{equation}\label{lema:mom_altern_proof_1}
\begin{split}
\frac{1}{\sqrt{k_n}}C_j^n(\widehat{\beta}_n) &= \frac{n}{k_n}\int_{\frac{(j-1)k_n}{n}}^{\frac{jk_n}{n}}(\beta_s-\widehat{\beta}_n)\sigma_s^2ds\\&~~+\frac{n}{k_n}\int_{\frac{(j-1)k_n}{n}}^{\frac{jk_n}{n}}X_s^ndY_s^c
+\frac{n}{k_n}\int_{\frac{(j-1)k_n}{n}}^{\frac{jk_n}{n}}Y_s^ndX_s^c-2\widehat{\beta}_n\frac{n}{k_n}\int_{\frac{(j-1)k_n}{n}}^{\frac{jk_n}{n}}X_s^ndX_s^c.\nonumber
\end{split}
\end{equation}
From here using the definition of the set $\mathcal{B}^n$ and applying the Burkholder-Davis-Gundy inequality, we have the result in \refeq{lema:mom_altern_1}. The results in \refeq{lema:mom_altern_2} and \refeq{lema:mom_altern_3} are shown in exactly the same way. Finally, the bounds on the probabilities in \refeq{lema:mom_altern_4} and \refeq{lema:mom_altern_5} follow from the fact that on $\mathcal{B}^n$, $\widehat{\beta}_n$ is bounded as well as an application of the Burkholder-Davis-Gundy inequality.
\qedd

\begin{lema}\label{lema:jumps}
Under Assumption SA, and with $n^{\tilde{\iota}}/k_n\rightarrow 0$ for some $\tilde{\iota}>0$, for any bounded random variable $b$ and $n$ sufficiently high, we have
\begin{equation}\label{lema:jumps_1}
\begin{split}
&\mathbb{E}_{\frac{(j-1)k_n}{n}}\left|\widehat{V}_j^{(n,1)}-V_j^{(n,1)}\right|^p+\mathbb{E}_{\frac{(j-1)k_n}{n}}\left|\widehat{V}_j^{(n,2)}(b)-V_j^{(n,2)}(b)\right|
^p\\&~~~~~~~~~~~~~~~~~~~~~\leq K\left(\frac{n^{p-1-(2p-r)\varpi}}{k_n^{p-1}}\bigvee n^{-p(2-r)\varpi}\right),~~~\textrm{for $p=1$, $p=2$ and $p=4$},
\end{split}
\end{equation}
\begin{equation}\label{lema:jumps_2}
\mathbb{E}_{\frac{(j-1)k_n}{n}}\left|\widehat{C}_j^n(b)-C_j(b)\right|^2\leq K\left(n^{1-(4-r)\varpi}\vee k_nn^{-2(2-r)\varpi}\right),
\end{equation}
\begin{equation}\label{lema:jumps_3}
\mathbb{P}\left(|\widehat{V}_j^{n}(\widehat{\beta}_n)-V_j^{n}(\widehat{\beta}_n)|1_{\{\mathcal{B}^n\}}>\epsilon\right)\leq Kn^{-(2-r)\varpi},~~\forall \epsilon>0.
\end{equation}
\end{lema}
\textbf{Proof of Lemma~\ref{lema:jumps}.}
In the proof we use the shorthand notation
\begin{equation}\label{lema:jumps_proof_1}
\mathcal{C}_i^n = \left\{|\Delta_i^nX|\leq\alpha\Delta_n^{\varpi},~~|\Delta_i^nY|\leq\alpha\Delta_n^{\varpi}\right\}.\nonumber
\end{equation}
We can decompose
\begin{equation}\label{lema:jumps_proof_2}
\begin{split}
\widehat{V}_j^{(n,1)}-V_j^{(n,1)} &= -\frac{n}{k_n}\sum_{i=(j-1)k_n+1}^{jk_n}(\Delta_i^nX^c)^21_{\{(\mathcal{C}_i^n)^c\}}+\frac{2n}{k_n}\sum_{i=(j-1)k_n+1}^{jk_n}\Delta_i^nX^c\Delta_i^nX^j
1_{\{\mathcal{C}_i^n\}}\\&~~~+\frac{n}{k_n}\sum_{i=(j-1)k_n+1}^{jk_n}(\Delta_i^nX^j)^21_{\{\mathcal{C}_i^n\}},\nonumber
\end{split}
\end{equation}
\begin{equation}\label{lema:jumps_proof_3}
\widehat{V}_j^{(n,2)}(b)-V_j^{(n,2)}(b) = \frac{n}{k_n}\sum_{i=(j-1)k_n+1}^{jk_n}\chi_i^{(n,1)}(b),~~\chi_i^{(n,1)}(b)=-(\Delta_i^nY^c-b\Delta_i^nX^c)^21_{\{(\mathcal{C}_i^n)^c\}},\nonumber
\end{equation}
\begin{equation}\label{lema:jumps_proof_4}
\chi_i^{(n,2)}(b)=(\Delta_i^nY^j-
b\Delta_i^nX^j)^21_{\{\mathcal{C}_i^n\}},~~\chi_i^{(n,3)}(b)=2(\Delta_i^nY^c-b\Delta_i^nX^c)(\Delta_i^nY^j-b\Delta_i^nX^j)1_{\{\mathcal{C}_i^n\}},\nonumber
\end{equation}
where we denoted $X^j = X-X^c$ and $Y^j = Y-Y^c$. We then have for $\forall p\geq 1$ and $\forall \iota>0$
\begin{equation}\label{lema:jumps_proof_5}
\left\{\begin{array}{l}\mathbb{E}_{i-1}^n\left[(\Delta_i^nX^c)^{2p}1_{\{(\mathcal{C}_i^n)^c\}}+|\chi_i^{(n,1)}(b)|^p\right]\leq Kn^{-p-1+r\varpi+\iota},\\
\mathbb{E}_{i-1}^n\left[(\Delta_i^nX^j)^{2p}1_{\{\mathcal{C}_i^n\}}+|\chi_i^{(n,2)}(b)|^p\right]\leq Kn^{-1-(2p-r)\varpi},\\ \mathbb{E}_{i-1}^n\left[|\Delta_i^nX^c\Delta_i^nX^j1_{\{\mathcal{C}_i^n\}}|^p+|\chi_i^{(n,3)}(b)|^p\right]\leq Kn^{-1-p/2-(p-r)\varpi+\iota}.\end{array}\right.
\end{equation}
Combining these results and using successive conditioning, we have the result in \refeq{lema:jumps_1}.

We next turn to \refeq{lema:jumps_2}. We have
\begin{equation}\label{lema:jumps_proof_6}
\widehat{C}_j^n(b)-C_j(b)=\frac{n}{\sqrt{k_n}}\sum_{i=(j-1)k_n+1}^{jk_n}\left(a_i+b_i+c_i+d_i\right),\nonumber
\end{equation}
\begin{equation}
\begin{split}
a_i &= -\Delta_i^nX^c(\Delta_i^nY^c-b\Delta_i^nX^c)1_{\{(\mathcal{C}_i^n)^c\}},~~b_i = \Delta_i^nX^c(\Delta_i^nY^j-b\Delta_i^nX^j)1_{\{\mathcal{C}_i^n\}},\\
c_i &= \Delta_i^nX^j(\Delta_i^nY^c-b\Delta_i^nX^c)1_{\{\mathcal{C}_i^n\}},~~d_i = \Delta_i^nX^j(\Delta_i^nY^j-b\Delta_i^nX^j)1_{\{\mathcal{C}_i^n\}}.\nonumber
\end{split}
\end{equation}
\begin{equation}\label{lema:jumps_proof_7}
\left\{\begin{array}{l}
\mathbb{E}_{i-1}^n|a_i|\leq Kn^{-2+r\varpi+\iota},~~\mathbb{E}_{i-1}^n|a_i|^2\leq Kn^{-3+r\varpi+\iota},\\
\mathbb{E}_{i-1}^n(|b_i|+|c_i|)\leq Kn^{-3/2-(1-r)\varpi+\iota},~~\mathbb{E}_{i-1}^n(|b_i|^2+|c_i|^2)\leq Kn^{-2-(2-r)\varpi+\iota},\\
~~\mathbb{E}_{i-1}^n|d_i|\leq Kn^{-1-(2-r)\varpi},~~\mathbb{E}_{i-1}^n|d_i|^2\leq Kn^{-1-(4-r)\varpi}.\end{array}\right.
\end{equation}
From here, using successive conditioning, we have the result in \refeq{lema:jumps_2}. We finally show \refeq{lema:jumps_3}. For some sufficiently big constant $\delta>0$ and sufficiently high $n$, taking into account the definition of the set $\mathcal{B}^n$, we have
\begin{equation}\label{lema:jumps_proof_8}
\begin{split}
&\mathbb{P}\left(|\widehat{V}_j^{n}(\widehat{\beta}_n)-V_j^{n}(\widehat{\beta}_n)|1_{\{\mathcal{B}^n\}}>\epsilon\right)\leq \mathbb{P}\left(|V_j^{(n,1)}|+|V_j^{(n,2)}(\widehat{\beta}_n)|>\delta\right)
+K\mathbb{E}\left|\widehat{V}_j^{(n,1)}-V_j^{(n,1)}\right|\\&~~~~~~~~+K\frac{n}{k_n}\sum_{k=1}^2\sum_{i=(j-1)k_n+1}^{jk_n}\mathbb{E}|\chi_i^{(n,k)}
(\overline{\beta}-\epsilon)|
+K\frac{n}{k_n}\sum_{k=1}^2\sum_{i=(j-1)k_n+1}^{jk_n}\mathbb{E}|\chi_i^{(n,k)}(\overline{\beta}+\epsilon)|\\&~~~~~~~~+K\frac{n}{k_n}\sum_{i=(j-1)k_n+1}^{jk_n}
\mathbb{E}\left[ \left(|\Delta_i^nX^c|+|\Delta_i^nY^c|\right)\left(|\Delta_i^nX^j|+|\Delta_i^nY^j|\right)1_{\{\mathcal{C}_i^n\}} \right],\nonumber
\end{split}
\end{equation}
and note that $\overline{\beta}$ coincides with $\beta_0$ under the null hypothesis and is a bounded positive random variable otherwise under Assumption SA. From here, applying the bounds in \refeq{lema:jumps_proof_5} above, as well as the bounds in \refeq{lema:mom_altern_4}-\refeq{lema:mom_altern_5} of Lemma~\ref{lema:mom_altern} and taking into account the rate of growth of $k_n$, we get the result in \refeq{lema:jumps_3}.
\qedd

\subsection{Proof of parts (a) of Theorems~\ref{thm:known_beta} and \ref{thm:unknown_beta}.}

We first prove the result for the statistic $T^n(\widehat{\beta}_n)$ with the result stated in the following lemma.

\begin{lema}\label{lema:cont}
Under Assumption SA with $\beta_t=\beta_0$ for $t\in[0,1]$ and further $k_n^{-1}n^{1/4}\rightarrow 0$ and $k_n^{-1}n\rightarrow \infty$, we have
\begin{equation}\label{lema:cont_1}
T^n(\widehat{\beta}_n)~\stackrel{\mathcal{L}}{\longrightarrow}~Z,\nonumber
\end{equation}
for $Z$ being a standard normal random variable.
\end{lema}
\textbf{Proof of Lemma~\ref{lema:cont}.}
We denote the sets
\begin{equation}\label{lema:cont_proof_1}
\mathcal{A}_{j}^{(n,1)} = \left\{\left| V_{j}^{(n,1)} - \frac{n}{k_n}\int_{\frac{(j-1)k_n}{n}}^{\frac{jk_n}{n}}\sigma_{s}^2ds \right|< \frac{1}{2} \frac{n}{k_n} \int_{\frac{(j-1)k_n}{n}}^{\frac{jk_n}{n}}\sigma_{s}^2ds\right\},\nonumber
\end{equation}
\begin{equation}\label{lema:cont_proof_2}
\mathcal{A}_{j}^{(n,2)}(b) = \left\{\left| V_{j}^{(n,2)}(b) - \frac{n}{k_n}\int_{\frac{(j-1)k_n}{n}}^{\frac{jk_n}{n}}\widetilde{\sigma}_{s}^2ds \right|< \frac{1}{2} \frac{n}{k_n} \int_{\frac{(j-1)k_n}{n}}^{\frac{jk_n}{n}}\widetilde{\sigma}_{s}^2ds\right\}.\nonumber
\end{equation}
We decompose
\begin{equation}\label{lema:cont_proof_3}
T_j^n(\widehat{\beta}_n)-T_j^n(\beta_0) = R_j^{(n,1)}+R_j^{(n,2)}+R_j^{(n,3)}+R_j^{(n,4)},\nonumber
\end{equation}
\begin{equation}\label{lema:cont_proof_4}
R^{(n,1)}_j = \frac{\left[(C_j^n(\widehat{\beta}_n)^2-V_j^n(\widehat{\beta}_n))-(C_j^n(\beta_0)^2-V_j^n(\beta_0))\right]\left(V_{j-1}^n(\beta_0)-V_{j-1}^n(\widehat{\beta}_n)\right)}
{\sigma^2_{\frac{(j-1)k_n}{n}}\widetilde{\sigma}^2_{\frac{(j-1)k_n}{n}}V_{j-1}^n(\widehat{\beta}_n) },\nonumber
\end{equation}
\begin{equation}\label{lema:cont_proof_5}
R^{(n,2)}_j = \frac{\left[(C_j^n(\widehat{\beta}_n)^2-V_j^n(\widehat{\beta}_n))-(C_j^n(\beta_0)^2-V_j^n(\beta_0))\right]\left(\sigma^2_{\frac{(j-1)k_n}{n}}
\widetilde{\sigma}^2_{\frac{(j-1)k_n}{n}}-V_{j-1}^n(\beta_0)\right)}{\sigma^2_{\frac{(j-1)k_n}{n}}\widetilde{\sigma}^2_{\frac{(j-1)k_n}{n}}V_{j-1}^n(\widehat{\beta}_n)},
\nonumber
\end{equation}
\begin{equation}\label{lema:cont_proof_6}
R^{(n,3)}_j = \frac{(C_j^n(\widehat{\beta}_n)^2-V_j^n(\widehat{\beta}_n))-(C_j^n(\beta_0)^2-V_j^n(\beta_0))}{\sigma^2_{\frac{(j-1)k_n}{n}}\widetilde{\sigma}^2_{\frac{(j-1)k_n}{n}}},
\nonumber
\end{equation}
\begin{equation}\label{lema:cont_proof_7}
R^{(n,4)}_j = \frac{(C_j^n(\beta_0)^2-V_j^n(\beta_0))\left(V_{j-1}^n(\beta_0)-V_{j-1}^n(\widehat{\beta}_n)\right)}{V_{j-1}^n(\beta_0)V_{j-1}^n(\widehat{\beta}_n)}.\nonumber
\end{equation}
We can further split
\begin{equation}\label{lema:cont_proof_8}
T_j^n(\beta_0) = T_j^{(n,1)}(\beta_0)+T_j^{(n,2)}(\beta_0),~~T_j^{(n,1)}(\beta_0) = \frac{C_j^n(\beta_0)^2-V_j^n(\beta_0)}{\sigma_{\frac{(j-1)k_n}{n}}^2\widetilde{\sigma}_{\frac{(j-1)k_n}{n}}^2},\nonumber
\end{equation}
\begin{equation}\label{lema:cont_proof_9}
T_j^{(n,2)}(\beta_0) =
\left(C_j^n(\beta_0)^2-V_j^n(\beta_0)\right)\left( \frac{1}{V_{j-1}^n(\beta_0)} -\frac{1}{\sigma_{\frac{(j-1)k_n}{n}}^2\widetilde{\sigma}_{\frac{(j-1)k_n}{n}}^2}\right).\nonumber
\end{equation}
We split the proof into several steps.

\noindent \textit{Step 1.} We prove $T^{(n,1)}(\beta_0)~\stackrel{\mathcal{L}}{\longrightarrow}~Z$. The result follows from an application of Theorem VIII.3.6 in \cite{JS}. In particular, using Lemmas~\ref{lema_mean}, \ref{lema_mean_2}, \ref{lema_var} and \ref{lema_powers}, we have
\begin{equation}\label{lema:cont_proof_10}
\left\{\begin{array}{l}\sqrt{\frac{k_n}{n}}\sum_{j=2}^{\lfloor n/k_n \rfloor} \mathbb{E}_{\frac{(j-1)k_n}{n}}\left(T_j^{(n,1)}(\beta_0)\right)~\stackrel{\mathbb{P}}{\longrightarrow}~0,\\
\frac{1}{2}\frac{k_n}{n}\sum_{j=2}^{\lfloor n/k_n \rfloor}\left\{ \mathbb{E}_{\frac{(j-1)k_n}{n}}\left(T_j^{(n,1)}(\beta_0)\right)^2-
\left(\mathbb{E}_{\frac{(j-1)k_n}{n}}\left(T_j^{(n,1)}(\beta_0)\right)\right)^2  \right\}~\stackrel{\mathbb{P}}{\longrightarrow}~1,\\
\sum_{j=2}^{\lfloor n/k_n \rfloor}\mathbb{P}\left( \left|\sqrt{\frac{k_n}{n}}T_j^{(n,1)}(\beta_0)\right|>\epsilon  \bigg|\mathcal{F}_{\frac{(j-1)k_n}{n}}\right)~\stackrel{\mathbb{P}}{\longrightarrow}~0,~~~\forall \epsilon>0.
\end{array}\right.\nonumber
\end{equation}

\noindent \textit{Step 2.} We prove $\sqrt{\frac{k_n}{n}}\sum_{j=2}^{\lfloor n/k_n\rfloor}\left(T_j^{(n,2)}(\beta_0)1_{\left\{\mathcal{A}_{j-1}^{(n,1)}~\cap~\mathcal{A}_{j-1}^{(n,2)}(\beta_0)\right\}}\right)~\stackrel{\mathbb{P}}{\longrightarrow}~0$.
Using Lemmas~\ref{lema_mean}, \ref{lema_mean_2} and \ref{lema_powers}, successive conditioning, the \Ito semimartingale assumption for $\sigma$ and $\widetilde{\sigma}$, and the definition of the sets $\mathcal{A}_j^{(n,1)}$ and $\mathcal{A}_j^{(n,2)}(\beta_0)$, we have 
\begin{equation}\label{lema:cont_proof_11}
\sum_{j=2}^{\lfloor n/k_n\rfloor}\mathbb{E}\left(T_j^{(n,2)}(\beta_0)1_{\left\{\mathcal{A}_{j-1}^n\right\}}\right)^2
\leq K\left(\frac{n}{k_n^2}\bigvee 1\right),\nonumber
\end{equation}
\begin{equation}\label{lema:cont_proof_13}
\left|\sum_{i,j:~i\neq j,~i\geq 2,~j\geq 2}\mathbb{E}\left(T_i^{(n,2)}(\beta_0)1_{\left\{\mathcal{A}_{i-1}^n\right\}}
T_j^{(n,2)}(\beta_0)1_{\left\{\mathcal{A}_{j-1}^n\right\}}  \right)\right|\leq K\left(\frac{n^{3/2}}{k_n^3}\bigvee1\right),\nonumber
\end{equation}
where we use the shorthand notation $\mathcal{A}_j^n = \mathcal{A}_j^{(n,1)}~\cap~\mathcal{A}_j^{(n,2)}(\beta_0)$. Combining the above three bounds, and taking into account the rate of growth condition for $k_n$, we establish the asymptotic negligibility result of this step.

\noindent \textit{Step 3.} We prove $\sqrt{\frac{k_n}{n}}\sum_{j=2}^{\lfloor n/k_n\rfloor}\left(R_j^{(n,1)}1_{\left\{\mathcal{A}_{j-1}^{(n,1)}~\cap~\mathcal{A}_{j-1}^{(n,2)}(\widehat{\beta}_n)\right\}}\right)~
\stackrel{\mathbb{P}}{\longrightarrow}~0$. First, we can decompose
\begin{equation}\label{lema:cont_proof_14}
C_j^n(\widehat{\beta}_n)^2-C_j^n(\beta_0)^2 = k_n(\widehat{\beta}_n-\beta_0)^2\left(V_j^{(n,1)}\right)^2-2\sqrt{k_n}(\widehat{\beta}_n-\beta_0)C_j^n(\beta_0)V_j^{(n,1)},
\end{equation}
\begin{equation}\label{lema:cont_proof_15}
V_j^n(\widehat{\beta}_n)-V_j^n(\beta_0) = \left(\widehat{\beta}_n-\beta_0\right)^2\left(V_j^{(n,1)}\right)^2-\frac{2}{\sqrt{k_n}}(\widehat{\beta}_n-\beta_0)V_j^{(n,1)}C_j^n(\beta_0).
\end{equation}
Next, using successive conditioning, Cauchy-Schwarz inequality, as well as the bounds derived in Lemma~\ref{lema_powers}, we get
\begin{equation}\label{lema:cont_proof_16}
\begin{split}
&\frac{k_n}{n}\sum_{j=2}^{\lfloor n/k_n\rfloor}\bigg[\mathbb{E}\left(|C_j^n(\beta_0)|V_j^{(n,1)}(V_{j-1}^{(n,1)})^2\right) + \mathbb{E}\left(|C_j^n(\beta_0)|V_j^{(n,1)}|C_{j-1}^n(\beta_0)|V_{j-1}^{(n,1)}\right)\\&~~~~~~~~~~~~~~~~~~~+ \mathbb{E}\left((V_j^{(n,1)})^2(V_{j-1}^{(n,1)})^2\right)+
\mathbb{E}\left((V_j^{(n,1)})^2|C_{j-1}^n(\beta_0)|V_{j-1}^{(n,1)}\right)\bigg]\leq K,\nonumber
\end{split}
\end{equation}
and the result of this step then follows from the $\sqrt{n}$ rate of convergence of $\widehat{\beta}_n$ to $\beta_0$ established in Lemma~\ref{lema:beta_hat} and the assumed conditions on $k_n$.

\noindent \textit{Step 4.} We prove $\sqrt{\frac{k_n}{n}}\sum_{j=2}^{\lfloor n/k_n\rfloor}\left(R_j^{(n,2)}1_{\left\{\mathcal{A}_{j-1}^{(n,1)}~\cap~\mathcal{A}_{j-1}^{(n,2)}(\widehat{\beta}_n)\right\}}\right)~
\stackrel{\mathbb{P}}{\longrightarrow}~0$. Using successive conditioning, the \Ito semimartingale assumption for $\sigma$ and $\widetilde{\sigma}$, as well as Lemma~\ref{lema_mean_2}, we get
\begin{equation}\label{lema:cont_proof_17}
\frac{k_n}{n}\sum_{j=2}^{\lfloor n/k_n\rfloor}\mathbb{E}\left\{\left[|C_j^n(\beta_0)|V_j^{(n,1)}+(V_j^{(n,1)})^2\right]\left|V_{j-1}^n(\beta_0) - \sigma^2_{\frac{(j-1)k_n}{n}}\widetilde{\sigma}^2_{\frac{(j-1)k_n}{n}}\right|\right\}\leq K\left(\sqrt{\frac{k_n}{n}}\bigvee\frac{1}{\sqrt{k_n}}\right),\nonumber
\end{equation}
and from here the result to be proved in this step follows because of the $\sqrt{n}$ rate of convergence of $\widehat{\beta}_n$ to $\beta_0$ established in Lemma~\ref{lema:beta_hat} and the decomposition in \refeq{lema:cont_proof_14}-\refeq{lema:cont_proof_15}.

\noindent \textit{Step 5.} We prove $\sqrt{\frac{k_n}{n}}\sum_{j=2}^{\lfloor n/k_n\rfloor}R_j^{(n,3)}~\stackrel{\mathbb{P}}{\longrightarrow}~0$. Given the CLT result in Lemma~\ref{lema:beta_hat} for $\widehat{\beta}_n$, and the decomposition in \refeq{lema:cont_proof_14}-\refeq{lema:cont_proof_15}, the result to be proved in this step will follow if we can show
\begin{equation}\label{lema:cont_proof_18}
\left(\frac{k_n}{n}\right)^{3/2}\sum_{j=2}^{\lfloor n/k_n\rfloor}\left(V_j^{(n,1)}\right)^2~\stackrel{\mathbb{P}}{\longrightarrow}~0,~~\frac{k_n}{n}\sum_{j=2}^{\lfloor n/k_n\rfloor}\left(C_j^n(\beta_0)V_j^{(n,1)}\right)~\stackrel{\mathbb{P}}{\longrightarrow}~0.\nonumber
\end{equation}
The first of this results follows trivially from the bound on the moments of $V_j^{(n,1)}$ derived in Lemma~\ref{lema_powers}. Next, application of Cauchy-Schwarz inequality, the proof of Lemma~\ref{lema_mean_2} and the bound of Lemma~\ref{lema_powers} for the $p$-th absolute moment of $C_j^n(\beta_0)$, yields
\[\frac{k_n}{n}\sum_{j=2}^{\lfloor n/k_n\rfloor}\left(C_j^n(\beta_0)\left(V_j^{(n,1)}-\sigma^2_{\frac{(j-1)k_n}{n}}\right)\right)~\stackrel{\mathbb{P}}{\longrightarrow}~0.\]
Finally, using the bound in \refeq{lema_mean_proof_3} in the proof of Lemma~\ref{lema_mean}, as well as successive conditioning and Cauchy-Schwarz and Burkholder-Davis-Gundy inequalities, we get
\begin{equation}\label{lema:cont_proof_19}
\mathbb{E}\left(\sum_{j=2}^{\lfloor n/k_n\rfloor}C_j^n(\beta_0)\right)^2\leq K\frac{n}{k_n},\nonumber
\end{equation}
and this implies the asymptotic negligibility of $\frac{k_n}{n}\sum_{j=2}^{\lfloor n/k_n\rfloor}C_j^n(\beta_0)$ and hence the result to be shown in this step.

\noindent \textit{Step 6.} We prove $\sqrt{\frac{k_n}{n}}\sum_{j=2}^{\lfloor n/k_n\rfloor}\left(R_j^{(n,4)}1_{\left\{\mathcal{A}_{j-1}^{(n,1)}~\cap~\mathcal{A}_{j-1}^{(n,2)}(\widehat{\beta}_n)~\cap~\mathcal{A}_{j-1}^{(n,2)}(\beta_0)\right\}}\right)~
\stackrel{\mathbb{P}}{\longrightarrow}~0$. Using successive conditioning, Cauchy-Schwarz inequality, as well as the bounds in Lemma~\ref{lema_powers}, we get
\begin{equation}\label{lema:cont_proof_20}
\frac{k_n}{n}\sum_{j=2}^{\lfloor n/k_n\rfloor}\mathbb{E}\left\{\left|(C_j^n(\beta_0))^2-V_j^{n}(\beta_0)\right|\left[(V_{j-1}^{(n,1)})^2 + |C_{j-1}^n(\beta_0)|V_{j-1}^{(n,1)}\right]\right\}\leq K,\nonumber
\end{equation}
and this implies the result to be shown in this step, given the $\sqrt{n}$ rate of convergence of $\widehat{\beta}_n$ and the decomposition in \refeq{lema:cont_proof_14}-\refeq{lema:cont_proof_15}.

\noindent \textit{Step 7.} We prove
\begin{equation}\label{lema:cont_proof_21}
\sqrt{\frac{k_n}{n}}\sum_{j=2}^{\lfloor n/k_n\rfloor}\left(\left(|T_j^{(n,2)}(\beta_0)|+|R_j^{(n,1)}|+|R_j^{(n,2)}|+|R_j^{(n,4)}|\right)
1_{\{\widetilde{A}_{j-1}^n\}}\right)~
\stackrel{\mathbb{P}}{\longrightarrow}~0,\nonumber
\end{equation}
for $\widetilde{\mathcal{A}}_j^n = \left\{\mathcal{A}_j^{(n,1)}~\cap~\mathcal{A}_j^{(n,2)}(\beta_0)~\cap~\mathcal{A}_j^{(n,2)}(\widehat{\beta}_n)\right\}^c$. In view of Lemma~\ref{lema:beta_hat}, it suffices to prove convergence on the set $\mathcal{B}^n$ for some sufficiently small positive numbers $\iota>0$ and $\delta>0$. We have
\begin{equation}\label{lema:cont_proof_22}
\begin{split}
&\mathbb{E}\left[\left|\sqrt{\frac{k_n}{n}}\sum_{j=2}^{\lfloor n/k_n\rfloor}\left(\left(|T_j^{(n,2)}(\beta_0)|+|R_j^{(n,1)}|+|R_j^{(n,2)}|+|R_j^{(n,4)}|\right)
1_{\{\widetilde{\mathcal{A}}_{j-1}^n\bigcap\mathcal{B}^n\}}\right)
\right|\bigwedge 1\right]\\&~~~\leq \sum_{j=2}^{\lfloor n/k_n\rfloor}\mathbb{P}\left(\widetilde{\mathcal{A}}_{j-1}^n
\bigcap\mathcal{B}^n\right)\leq K\frac{n}{k_n}\frac{1}{k_n^{p/2}},~~\forall p\geq 1,\nonumber
\end{split}
\end{equation}
where for the last inequality, we applied Lemmas~\ref{lema_powers} and \ref{lema:mom_altern} and the boundedness from below of the processes $|\sigma|$ and $|\widetilde{\sigma}|$. This bound implies the asymptotic negligibility to be proved in this step.
\qedd

\medskip
\noindent \textbf{Proof of parts (a) of Theorems~\ref{thm:known_beta} and \ref{thm:unknown_beta} continued.}
What remains to be shown is that the difference $\widehat{T}^n(\widehat{\beta}_n)-T^n(\widehat{\beta}_n)$ is asymptotically negligible. Recalling the definition of the set $\mathcal{B}^n$, we first note that $\mathbb{P}((\mathcal{B}^n)^c)\rightarrow~0$, therefore it suffices to focus on the set $\mathcal{B}^n$ only.

We decompose $\widehat{T}^n_j(\widehat{\beta}_n)-T^n_j(\widehat{\beta}_n) = \xi_j^{(n,1)}+\xi_j^{(n,2)}+\xi_j^{(n,3)}+\xi_j^{(n,4)}$, where
\begin{equation}\label{proof_a_1}
\xi_j^{(n,1)} = \frac{ \widehat{C}_j^n(\widehat{\beta}_n)^2-C_j^n(\widehat{\beta}_n)^2-\widehat{V}_{j}^n(\widehat{\beta}_n)+V_{j}^n(\widehat{\beta}_n) }{\widehat{V}_{j-1}^n(\widehat{\beta}_n)},\nonumber
\end{equation}
\begin{equation}\label{proof_a_2}
\xi_j^{(n,2)} = \frac{ C_j^n(\widehat{\beta}_n)^2-C_j^n(\beta_0)^2-V_{j}^n(\widehat{\beta}_n)+V_{j}^n(\beta_0) }{V_{j-1}^n(\widehat{\beta}_n)\widehat{V}_{j-1}^n(\widehat{\beta}_n)}\left(V_{j-1}^n(\widehat{\beta}_n)-\widehat{V}_{j-1}^n(\widehat{\beta}_n)\right),\nonumber
\end{equation}
\begin{equation}\label{proof_a_3}
\xi_j^{(n,3)} = \frac{C_j^n(\beta_0)^2-V_{j}^n(\beta_0) }{V_{j-1}^n(\beta_0)\widehat{V}_{j-1}^n(\beta_0)}\left(V_{j-1}^n(\beta_0)-\widehat{V}_{j-1}^n(\beta_0)\right),
\nonumber
\end{equation}
\begin{equation}\label{proof_a_4}
\xi_j^{(n,4)} = (C_j^n(\beta_0)^2-V_{j}^n(\beta_0))\left(\frac{1}{\widehat{V}_{j-1}^n(\widehat{\beta}_n)}-\frac{1}{\widehat{V}_{j-1}^n(\beta_0)}-\frac{1}{V_{j-1}^n(\widehat{\beta}_n)}
+\frac{1}{V_{j-1}^n(\beta_0)}\right).\nonumber
\end{equation}
The proof consists of several steps and we will henceforth denote with $\epsilon$ some sufficiently small positive constant.

\noindent \textit{Step 1.} We prove $1_{\{\mathcal{B}^n\}}\sqrt{\frac{k_n}{n}}\sum_{j=2}^{\lfloor n/k_n \rfloor}\xi_j^{(n,1)}1_{\{|\widehat{V}_{j-1}^n(\widehat{\beta}_n)|>\epsilon\}}~\stackrel{\mathbb{P}}{\longrightarrow}~0$,
whenever $n^{1/2-(2-r)\varpi}\rightarrow 0$ and $\frac{n^{2-(4-r)\varpi}}{k_n}\rightarrow 0$. This follows directly from applying the algebraic identity $x^2-y^2=(x-y)^2+2y(x-y)$ for any real $x$ and $y$, the bound on $\widehat{\beta}_n-\beta_0$ on the set $\mathcal{B}^n$, as well as Lemmas~\ref{lema_powers} and \ref{lema:jumps}.

\noindent \textit{Step 2.} We prove $1_{\{\mathcal{B}^n\}}\sqrt{\frac{k_n}{n}}\sum_{j=2}^{\lfloor n/k_n \rfloor}\xi_j^{(n,2)}1_{\{|\widehat{V}_{j-1}^n(\widehat{\beta}_n)|>\epsilon,~V_{j-1}^n(\widehat{\beta}_n)>\epsilon\}}~\stackrel{\mathbb{P}}{\longrightarrow}~0$. We make use of the decomposition in \refeq{lema:cont_proof_14} and \refeq{lema:cont_proof_15}, the bounds in Lemmas~\ref{lema_powers} and \ref{lema:jumps}, the fact that $n^{1/2-\iota}(\widehat{\beta}_n-\beta_0)~\stackrel{\mathbb{P}}{\longrightarrow}~0$ for arbitrary small $\iota>0$ (by Lemma~\ref{lema:beta_hat}).

\noindent \textit{Step 3.} We prove $1_{\{\mathcal{B}^n\}}\sqrt{\frac{k_n}{n}}\sum_{j=2}^{\lfloor n/k_n \rfloor}\xi_j^{(n,3)}1_{\{|\widehat{V}_{j-1}^n(\beta_0)|>\epsilon,~|V_{j-1}^n(\beta_0)|>\epsilon\}}~\stackrel{\mathbb{P}}{\longrightarrow}~0$, provided $\frac{n^{1-\frac{(8-r)\varpi}{3}}}{k_n}\rightarrow 0$.
This result follows from showing convergence in $L^2$-norm, upon applying successive conditioning, using Cauchy-Schwarz inequality and the bounds in Lemmas~\ref{lema_mean}, \ref{lema_var}, \ref{lema:mom_altern} and \ref{lema:jumps}.

\noindent \textit{Step 4.} We prove $1_{\{\mathcal{B}^n\}}\sqrt{\frac{k_n}{n}}\sum_{j=2}^{\lfloor n/k_n \rfloor}\xi_j^{(n,4)}1_{\{|\widehat{V}_{j-1}^n(\beta_0)|>\epsilon,~|V_{j-1}^n(\beta_0)|>\epsilon, ~|\widehat{V}_{j-1}^n(\widehat{\beta}_n)|>\epsilon,~|V_{j-1}^n(\widehat{\beta}_n)|>\epsilon\}}~\stackrel{\mathbb{P}}{\longrightarrow}~0$, provided
$\frac{n^{\frac{1}{3}-\frac{2}{3}(4-r)\varpi+\iota}}{k_n}\rightarrow 0$ for some arbitrary small $\iota>0$.

Making use of an analogous decomposition of the difference $\widehat{V}_{j}^n(\widehat{\beta}_n) - \widehat{V}_{j}^n(\beta_0)$ as in \refeq{lema:cont_proof_15}, we can bound
\begin{equation}\label{proof_a_5}
\begin{split}
&\left|\frac{1}{\widehat{V}_{j-1}^n(\widehat{\beta}_n)}-\frac{1}{\widehat{V}_{j-1}^n(\beta_0)}\right|1_{\{|\widehat{V}_{j-1}^n(\beta_0)|>\epsilon, ~|\widehat{V}_{j-1}^n(\widehat{\beta}_n)|>\epsilon\}}
\\&~~~~~~~~~~~~~~~~~~~~~~~~~~~~~~~~~~~~~~~\leq \frac{K}{\sqrt{k_n}}|\widehat{\beta}_n-\beta_0|\widehat{V}_j^{(n,1)}|\widehat{C}_j^n(\beta_0)|+K|\widehat{\beta}_n-\beta_0|^2(\widehat{V}_j^{(n,1)})^2.\nonumber
\end{split}
\end{equation}
Similar analysis can be made for the term involving $1/V_{j}^n(\widehat{\beta}_n) - 1/V_{j}^n(\beta_0)$. From here using the convergence result for $\widehat{\beta}_n$ in Lemma~\ref{lema:beta_hat}, and the results in Lemmas~\ref{lema_var}, \ref{lema_powers} and \ref{lema:jumps}, we get the result to be proved in this step.

\noindent \textit{Step 5.} We prove
\begin{equation}\label{proof_a_6}
\begin{split}
&1_{\{\mathcal{B}^n\}}\sqrt{\frac{k_n}{n}}\sum_{j=2}^{\lfloor n/k_n \rfloor}\bigg[(|\xi_j^{(n,1)}|+|\xi_j^{(n,2)}|+|\xi_j^{(n,3)}|+|\xi_j^{(n,4)}|)\\&~~~~~~~~~~~~~~~~~~~~~~~~~\times1_{\{|\widehat{V}_{j-1}^n(\beta_0)|<\epsilon~\cup~
|V_{j-1}^n(\beta_0)|<\epsilon~\cup~|\widehat{V}_{j-1}^n(\widehat{\beta}_n)|<\epsilon~\cup~|V_{j-1}^n(\widehat{\beta}_n)|<\epsilon\}}\bigg]~\stackrel{\mathbb{P}}
{\longrightarrow}~0,\nonumber
\end{split}
\end{equation}
provided $\frac{n^{1-(2-r)\varpi}}{k_n}\rightarrow 0$. Using the definition of the set $\mathcal{B}^n$, the bounds in Lemmas~\ref{lema_powers}, \ref{lema:mom_altern} and \ref{lema:jumps}, we get the result of this step.
\qedd

\subsection{Proof of parts (b) of Theorems~\ref{thm:known_beta} and \ref{thm:unknown_beta}.}
We can decompose $\frac{1}{k_n}\widehat{T}_j^n(\widehat{\beta}_n) = \overline{\xi}_j^{(n,1)}+\overline{\xi}_j^{(n,2)}+\overline{\xi}_j^{(n,3)}$, where
\begin{equation}\label{proof_b_1}
\overline{\xi}_j^{(n,1)} = \frac{\left(\frac{n}{k_n}\int_{\frac{(j-1)k_n}{n}}^{\frac{jk_n}{n}}(\beta_{s}-\overline{\beta})\sigma_{s}^2ds\right)^2  }{ \frac{n}{k_n}\int_{\frac{(j-1)k_n}{n}}^{\frac{jk_n}{n}}\sigma_{s}^2ds\frac{n}{k_n}\int_{\frac{(j-1)k_n}{n}}^{\frac{jk_n}{n}}\left( (\beta_{s}-\overline{\beta})^2\sigma_{s}^2+\widetilde{\sigma}^2_{s} \right)ds  },\nonumber
\end{equation}
\begin{equation}\label{proof_b_2}
\overline{\xi}_j^{(n,2)} = \frac{\frac{1}{k_n}\left(\widehat{C}_j^n(\widehat{\beta}_n)^2-\widehat{V}_j^n(\widehat{\beta}_n) \right) -\left(\frac{n}{k_n}\int_{\frac{(j-1)k_n}{n}}^{\frac{jk_n}{n}}(\beta_{s}-\overline{\beta})\sigma_{s}^2ds\right)^2  }{ \frac{n}{k_n}\int_{\frac{(j-1)k_n}{n}}^{\frac{jk_n}{n}}\sigma_{s}^2ds\frac{n}{k_n}\int_{\frac{(j-1)k_n}{n}}^{\frac{jk_n}{n}}\left( (\beta_{s}-\overline{\beta})^2\sigma_{s}^2+\widetilde{\sigma}^2_{s} \right)ds  },\nonumber
\end{equation}
\begin{equation}\label{proof_b_3}
\overline{\xi}_j^{(n,3)} = \frac{1}{k_n}\frac{\widehat{C}_j^n(\widehat{\beta}_n)^2-\widehat{V}_j^n(\widehat{\beta}_n)}{\widehat{V}_{j-1}^n(\widehat{\beta}_n)}-
\frac{1}{k_n}\frac{\widehat{C}_j^n(\widehat{\beta}_n)^2-\widehat{V}_j^n(\widehat{\beta}_n)}{ \frac{n}{k_n}\int_{\frac{(j-1)k_n}{n}}^{\frac{jk_n}{n}}\sigma_{s}^2ds\frac{n}{k_n}\int_{\frac{(j-1)k_n}{n}}^{\frac{jk_n}{n}}\left( (\beta_{s}-\overline{\beta})^2\sigma_{s}^2+\widetilde{\sigma}^2_{s} \right)ds  }.\nonumber
\end{equation}
The proof consists of the following steps in which we denote with $\epsilon$ some sufficiently small positive constant.

\noindent \textit{Step 1.} We have $\frac{k_n}{n}\sum_{j=2}^{\lfloor n/k_n \rfloor}\overline{\xi}_j^{(n,1)}~\stackrel{a.s.}{\longrightarrow}~\int_0^1\frac{(\beta_{s}-\overline{\beta})^2\sigma_{s}^2}
{\left((\beta_{s}-\overline{\beta})^2
\sigma_{s}^2+\widetilde{\sigma}^2_{s}\right)}ds$. This follows from convergence of Riemann sums and the fact that the processes $\beta$, $\sigma$ and $\widetilde{\sigma}$ have \cadlag paths (and hence are Riemann integrable).

\noindent \textit{Step 2.} We prove $\frac{k_n}{n}\sum_{j=2}^{\lfloor n/k_n \rfloor}\overline{\xi}_j^{(n,2)}~\stackrel{\mathbb{P}}{\longrightarrow}~0$, provided $\frac{n^{1-(4-r)\varpi}}{k_n}\rightarrow 0$. First, it suffices to focus attention to the set $\mathcal{B}^n$ because of the CLT for $\widehat{\beta}_n$ in Lemma~\ref{lema:beta_hat} and we do so. Then using the bounds in Lemma~\ref{lema:mom_altern}, we have $\mathbb{E}|\overline{\xi}_j^{(n,2)}|\leq K\left(\frac{1}{\sqrt{k_n}}\bigvee \frac{n^{\frac{1-(4-r)\varpi}{2}}}{\sqrt{k_n}}\bigvee\frac{n^{1-(4-r)\varpi}}{k_n}\bigvee n^{-(2-r)\varpi}\right)$, and this implies the result to be shown in this step.

\noindent \textit{Step 3.} We prove $\frac{k_n}{n}\sum_{j=2}^{\lfloor n/k_n \rfloor}\overline{\xi}_j^{(n,3)}1_{\{|\widehat{V}_{j-1}^n(\widehat{\beta}_n)|>\epsilon\}}~\stackrel{\mathbb{P}}{\longrightarrow}~0$, provided $\frac{n^{1-(4-r)\varpi}}{k_n}\rightarrow 0$.
First, we have the algebraic inequality:
\begin{equation}\label{proof_b_4}
\begin{split}
&\left|\frac{1}{k_n}\widehat{C}_j^n(\widehat{\beta}_n)^2 - \left(\frac{n}{k_n}\int_{\frac{(j-1)k_n}{n}}^{\frac{jk_n}{n}}(\beta_{s}-\overline{\beta})\sigma_{s}^2ds\right)^2\right|
\widehat{V}_{j-1}^n(\widehat{\beta}_n)1_{\{\mathcal{B}^n\}}\\&~~~~~~~~~~~~~~~~~~~\leq K
\left|\frac{1}{k_n}\widehat{C}_j^n(\widehat{\beta}_n)^2 - \left(\frac{n}{k_n}\int_{\frac{(j-1)k_n}{n}}^{\frac{jk_n}{n}}(\beta_{s}-\overline{\beta})\sigma_{s}^2ds\right)^2\right|
\widehat{V}_{j-1}^{(n,1)}(\widehat{V}_{j-1}^{(n,1)}+\widehat{V}_{j-1}^{(n,2)}(0)).\nonumber
\end{split}
\end{equation}
From here, using successive conditioning and Lemmas~\ref{lema:mom_altern} and \ref{lema:jumps}, we get
\begin{equation}\label{proof_b_5}
\mathbb{E}\left(\left|\frac{1}{k_n}\widehat{C}_j^n(\widehat{\beta}_n)^2 - \left(\frac{n}{k_n}\int_{\frac{(j-1)k_n}{n}}^{\frac{jk_n}{n}}(\beta_{s}-\overline{\beta})\sigma_{s}^2ds\right)^2\right|
\widehat{V}_{j-1}^n(\widehat{\beta}_n)1_{\{\mathcal{B}^n\}}\right)\leq K\eta_n,\nonumber
\end{equation}
where we use the shorthand $\eta_n = \left(\frac{1}{\sqrt{k_n}}\bigvee \frac{n^{\frac{1-(4-r)\varpi}{2}}}{\sqrt{k_n}}\bigvee\frac{n^{1-(4-r)\varpi}}{k_n}\bigvee n^{-(2-r)\varpi}\right)$. Similar analysis leads to
\begin{equation}\label{proof_b_6}
\mathbb{E}\left(\left|\widehat{V}_j^n(\widehat{\beta}_n) - \frac{n}{k_n}\int_{\frac{(j-1)k_n}{n}}^{\frac{jk_n}{n}}\sigma_{s}^2ds\frac{n}{k_n}\int_{\frac{(j-1)k_n}{n}}^{\frac{jk_n}{n}}\left( (\beta_{s}-\overline{\beta})^2\sigma_{s}^2+\widetilde{\sigma}^2_{s} \right)ds\right|
\widehat{V}_{j-1}^n(\widehat{\beta}_n)1_{\{\mathcal{B}^n\}}\right)\leq K\eta_n.\nonumber
\end{equation}
Combining the above two results, it is hence sufficient to prove the result of this step in which $\frac{1}{k_n}\left(\widehat{C}_j^n(\widehat{\beta}_n)^2-\widehat{V}_j^n(\widehat{\beta}_n)\right)$ in $\overline{\xi}_j^{(n,3)}$ is replaced with the term 
\[\left(\frac{n}{k_n}\int_{\frac{(j-1)k_n}{n}}^{\frac{jk_n}{n}}(\beta_{s}-\overline{\beta})\sigma_{s}^2ds\right)^2-\frac{1}{k_n}\frac{n}{k_n}
\int_{\frac{(j-1)k_n}{n}}^{\frac{jk_n}{n}}\sigma_{s}^2ds\frac{n}{k_n}\int_{\frac{(j-1)k_n}{n}}^{\frac{jk_n}{n}}\left( (\beta_{s}-\overline{\beta})^2\sigma_{s}^2+\widetilde{\sigma}^2_{s} \right)ds.\] The last term is bounded and hence the result follows by an application of Lemmas~\ref{lema:mom_altern} and \ref{lema:jumps}.

\noindent \textit{Step 4.} We prove $\frac{k_n}{n}\sum_{j=2}^{\lfloor n/k_n \rfloor}\overline{\xi}_j^{(n,3)}1_{\{|\widehat{V}_{j-1}^n(\widehat{\beta}_n)|\leq\epsilon\}}~\stackrel{\mathbb{P}}{\longrightarrow}~0$ provided $\frac{n^{1-(2-r)\varpi}}{k_n}\rightarrow 0$. We will be done if we can show $\mathbb{P}\left(|\widehat{V}_{j-1}^n(\widehat{\beta}_n)|\leq\epsilon~\cap~\mathcal{B}^n\right)\leq K\left(\frac{k_n}{n}\right)\eta_n$ for some deterministic sequence $\eta_n\rightarrow 0$. This follows from an application of the bounds in Lemmas~\ref{lema:mom_altern} and \ref{lema:jumps}, as well as the fact that on $\mathcal{B}^n$ $\widehat{\beta}_n$ is bounded.
\qedd

\subsection{Proof of Theorem~\ref{thm:la_sup}.}

For a general process $\beta_t$ we define $Z_t^n=Y_t-\int_0^t\beta_{\lfloor sn\rfloor/n}dX_s$ and we split $\widehat C_j^n(\beta)=\widehat C_j^{n,X}(\beta)+\widehat C_j^{n,Z}(\beta)$ with
\begin{equation}
\left\{\begin{array}{l}\widehat C_j^{n,X}(\beta)= \frac{n}{\sqrt{k_n}}\sum_{i=(j-1)k_n+1}^{jk_n}s_i^2(\beta_{(i-1)/n}-\beta),~~s_i^2=(\Delta_i^nX)^2, \\ \widehat C_j^{n,Z}(\beta)= \frac{n}{\sqrt{k_n}}\sum_{i=(j-1)k_n+1}^{jk_n}\Delta_i^nX\Delta_i^nZ^n. \end{array}\right.
\end{equation}
We further set
\[ \widehat V_j^{n,X}(\beta)=\Big(\frac{n}{k_n}\sum_{i=(j-1)k_n+1}^{jk_n}s_i^2 \Big)\Big( \frac{n}{k_n}\sum_{i=(j-1)k_n+1}^{jk_n}s_i^2(\beta_{(i-1)/n}
 -\beta)^2\Big),
 \]
and
\[ \widehat T^{n,X}(\beta)=\frac{1}{\sqrt{2}} \sqrt{\frac{k_n}{n}}\sum_{j=2}^{\lfloor\frac{n}{k_n}\rfloor}\frac{\widehat C_j^{n,X}(\beta)^2- \widehat{V}_j^{n,X}(\beta)}{\widehat V_{j-1}^n(\beta)}.
\]
We start with the analysis of $\widehat T^{n,X}(\beta)$. First, by the H\"older property, we have
\begin{align*}
\widehat C_j^{n,X}(\beta)^2 &=\frac{n^2}{k_n}\Big(\sum_{i=(j-1)k_n+1}^{jk_n}s_i^2\Big)
 \Big(\sum_{i=(j-1)k_n+1}^{jk_n}s_i^2(\beta_{(i-1)/n}
 -\beta)^2 - \sum_{i=(j-1)k_n+1}^{jk_n}s_i^2(\beta_{(i-1)/n}
 -\bar\beta_j)^2\Big)\\
&\ge \frac{n^2}{k_n} \Big(\sum_{i=(j-1)k_n+1}^{jk_n}s_i^2\Big)\Big(\sum_{i=(j-1)k_n+1}^{jk_n}s_i^2(\beta_{(i-1)/n}
 -\beta)^2 - R^2(k_n/n)^{2\alpha}\sum_{i=(j-1)k_n+1}^{jk_n}s_i^2\Big),
\end{align*}
where we denoted $\bar \beta_j=\sum_{i=(j-1)k_n+1}^{jk_n}s_i^2\beta_{(i-1)/n}/\sum_{i=(j-1)k_n+1}^{jk_n}s_i^2$. From here, by straightforward expectation and variance bounds, and taking into account the assumed rate at which $k_n$ grows asymptotically, we find
\begin{align*}
\widehat C_j^{n,X}(\beta)^2-\widehat V_j^{n,X}(\beta)&\ge\frac{n^2}{k_n}(1-O_{L^2}(k_n^{-1/2}))\int_{(j-1)k_n/n}^{jk_n/n}\sigma_t^2dt\\
&\quad \times\int_{(j-1)k_n/n}^{jk_n/n}\sigma_t^2\Big((\beta_t-\beta)^2(1- O_{L^2}(k_n^{-1/2}))
-R^2(k_n/n)^{2\alpha}(1+O_{L^2}(k_n^{-1/2}))\Big)\,dt.
\end{align*}
It is easy to see that Lemma~\ref{lema:mom_altern} continues to hold with $Y_t$ replaced by  $\int_0^t\beta_{\lfloor sn\rfloor/n} dX_s$. Therefore, for any $\delta\in(0,1)$ there is a $K>0$ such that
\begin{align*}
 &\mathbb{P}\Big(\bigcup_{j=1,\ldots,\lfloor n/k_n \rfloor}D_j\Big)
\le K k_n^{-1}\text{ with }\\
 &D_j:=\Big\{  \frac{n}{k_n}\int_{(j-1)k_n/n}^{jk_n/n}\sigma_t^2dt\frac{n}{k_n}\int_{(j-1)k_n/n}^{jk_n/n}((\beta_t-\beta)^2\sigma_t^2+\widetilde{\sigma}_t^2)\,dt \big/\widehat V_j^n(\beta)\notin [1-\delta,1+\delta] \Big\}.
\end{align*}
Hence, we can asymptotically work on the event $\bigcap_{j=1,\ldots,\lfloor n/k_n \rfloor}D_j^c$ and conclude that $\widehat T^{n,X}(\beta)$ is bounded from below by
\[ \sqrt{\frac{nk_n}{2}}\sum_{j=2}^{\lfloor\frac{n}{k_n}\rfloor}
\frac{(1-\delta) \int_{(j-1)k_n/n}^{jk_n/n}\sigma_t^2(\beta_t-\beta)^2dt
-(1+\delta)R^2(k_n/n)^{2\alpha}\int_{(j-1)k_n/n}^{jk_n/n}\sigma_t^2dt}
{\frac{n}{k_n}\int_{(j-1)k_n/n}^{jk_n/n}((\beta_t-\beta)^2\sigma_t^2+\widetilde{\sigma}_t^2)\,dt} - O_P(1).
\]
The latter is, using Riemann sum approximations and the separation $\Gamma r_n$ under the alternative, of  order
\[ \sqrt{\frac{nk_n}{2}}\Big((1-\delta)\Gamma^2r_n^2
-(1+\delta)R^2(k_n/n)^{2\alpha}\Big)\int_{0}^1\frac{\sigma_t^2}{(\beta_t-\beta)^2\sigma_t^2+\widetilde{\sigma}_t^2}dt- O_P(1),
\]
where we also made use of the fact that $\sigma^2$ and $\widetilde{\sigma}^2$ are \Ito semimartingales as well as the assumed growth condition for $k_n$. Consequently, for  $\Gamma= \sqrt{(\sqrt{2}K+1)\frac{1+\delta}{1-\delta}} R$ and the choices of $k_n$ and $r_n$ we have
\[  \widehat T^{n,X}(\beta)\ge KR^2\int_{0}^1\frac{\sigma_t^2}{(\beta_t-\beta)^2\sigma_t^2+\widetilde{\sigma}_t^2}dt- O_P(1)\]
uniformly over the alternative and over $K>0$, $n\ge 1$.

The same arguments as for the proof of Theorem 1(a), applied to the pair $\left(X_t,Y_t-\int_0^t\beta_sdX_s\right)$, the fact that $d\langle X^c,(Z^n)^c\rangle_t=(\beta_t-\beta_{\lfloor tn\rfloor/n})\sigma_t^2dt$ is asymptotically negligible due to $\sup_{t\in[0,1]}|(\beta_t-\beta_{\lfloor tn\rfloor/n})\sigma_t^2|=O_P(n^{-\alpha})$ for $\beta_t\in C^\alpha(R)$, together with the bounds derived above and a Cauchy-Schwarz bound for the cross term yield the following uniform result
\[ \widehat T^{n}(\beta)\ge K R^2\Big(\int_0^1\frac{\sigma_t^2}{(\beta_t-\beta)\sigma_t^2+\widetilde{\sigma}_t^2}dt\Big)^2- O_P(\sqrt{K}).\]
The right-hand side converges to $+\infty$ in probability as $K\to+\infty$. Hence, we choose $K$ and thus $\Gamma$ so large that $P_{\beta_t}(\widehat T^{n}(\beta)\ge c_{\gamma/2})>1-\gamma/2$ on $H_{1,\alpha}(\Gamma r_n)$ holds, which implies the result.\qedd

\subsection{Proof of Theorem~\ref{thm:la_inf}.}
For any sign sequence $\eps=(\eps_j)\in\{-1,+1\}^{J_n}$ with $J_n=n/k_n$, we define
\[ \beta_\eps(t)=\beta+\sum_{j=0}^{J_n-1} \eps_j J_n^{-\alpha} K(J_nt-j) ,\]
where $K$ is a kernel of support $[0,1]$ with $\alpha$-H\"older constant smaller than $R$. Then the functions $J_n^{-\alpha} K(J_nt-j)$, $j=0,\ldots,J_n-1$, have disjoint support on $[0,1]$ and lie in $C^\alpha(R)$. Consequently, also $\beta_\eps(t)$ is in the H\"older ball $C^\alpha(R)$. We obtain further $\|\beta_\eps(t)-\beta\|_{L^2}=J_n^{-\alpha}\|K\|_{L^2}$ and thus $\beta_\eps(t)\in H^1(cJ_n^{-\alpha}\|K\|_{L^2})$ holds for some $c>0$.

Under the alternatives we work with some fixed (deterministic and positive) $\sigma^2_t$ and $\widetilde{\sigma}^2_t$. For the hypothesis $H_0:\beta_t=\beta$ we
set
\[ \rho_i=\frac{\int_{(i-1)/n}^{i/n} \sigma_t^2\abs{\beta_\eps(t)-\beta}\,dt}{\sqrt{\int_{(i-1)/n}^{i/n} \sigma_t^2dt}\sqrt{\int_{(i-1)/n}^{i/n} (\widetilde{\sigma}_t^2+\sigma_t^2(\beta_\eps(t)-\beta)^2)\,dt}},
\]
which is independent of $\eps$ because $\abs{\beta_\eps(t)-\beta}$ does not depend on the sign. Under $H_0$ we then consider volatilities $\sigma_{0,t}^2$ and $\widetilde{\sigma}_{0,t}^2$, depending on $n$, such that $\int_{(i-1)/n}^{i/n}\sigma_{0,t}^2dt=(1-\rho_i^2)\int_{(i-1)/n}^{i/n}\sigma_t^2dt$ and with $\widetilde{\sigma}_{0,t}^2$ defined in an analogous way. Note that by construction and by H\"older continuity of $\beta_\eps(t)$ we have $\rho_i=O(J_n^{-\alpha})$ and $\abs{\rho_i-\rho_{i-1}}=O(n^{-\alpha})$ (recall our assumption for $\sigma$ and $\widetilde{\sigma}$) so that for each $n$ we can even find a smooth version of $\sigma_{0}^2$ and $\widetilde{\sigma}_{0}^2$. This minimal change simplifies the ensuing likelihood considerations drastically because it guarantees that the empirical covariances are sufficient statistics for these sets of parameters.

We bound the minimax testing error by the average error over $\beta_\eps$ using the likelihood to change the measure and the Cauchy-Schwarz inequality in combination with $E_\beta[\psi_n^2]\le 1$, $E_\beta[\frac{dP_{\beta_\eps}}{dP_\beta}]=1$:
\begin{align*}
P_\beta(\psi_n=1)+\sup_{\beta_t\in H_{1,\alpha}(cJ_n^{-\alpha}\|K\|_{L^2})}P_{\beta_t}(\psi_n=0) & \ge 2^{-J_n}\sum_{\eps\in\{-1,+1\}^{J_n}}\big( P_\beta(\psi_n=1)+P_{\beta_\eps(t)}(\psi_n=0)\big)\\
& =1-E_\beta\Big[\psi_n\Big(2^{-J_n}\sum_{\eps\in\{-1,+1\}^{J_n}}\Big(\frac{dP_{\beta_\eps}}{dP_\beta}-1\Big)\Big)\Big]\\
& \ge 1-\Big(E_\beta\Big[\Big(2^{-J_n}\sum_{\eps\in\{-1,+1\}^{J_n}}\frac{dP_{\beta_\eps}}{dP_\beta}\Big)^2\Big]-1\Big)^{1/2}.
\end{align*}
Since the transformed increments $\Delta_i^n(X,Y-\beta X)=(\Delta_i^nX,\Delta_i^nY-\beta\Delta_i^nX)$ are independent under all $P_{\beta_t}$, the likelihood factorizes over the $J_n$ blocks:
\[ \frac{dP_{\beta_\eps}}{dP_\beta}(\Delta_i^n(X,Y-\beta X)_{1\le i\le n})=\prod_{j=0}^{J_n-1} \frac{p_{j,\eps_j}(\Delta_{jk_n+i}^n(X,Y-\beta X))_{1\le i\le k_n}}{p_{j,0}(\Delta_{jk_n+i}^n(X,Y-\beta X))_{1\le i\le k_n}}\]
with density functions $p_{j,1},p_{j,-1},p_{j,0}$ on $\R^{2k_n}$ of the transformed increments on block $j$. This factorization permits a significant simplification, using invariance with respect to bi-measurable transformations:
\[ E_\beta\Big[\Big(2^{-J_n}\sum_{\eps\in\{-1,+1\}^{J_n}}\frac{dP_{\beta_\eps}}{dP_\beta}\Big)^2\Big]= 2^{-2J_n}\sum_{\eps,\eps'\in\{-1,+1\}^{J_n}}\prod_{j=0}^{J_n-1}\int \frac{p_{j,\eps_j}p_{j,\eps'_j}}{p_{j,0}^2}p_{j,0}
= \prod_{j=0}^{J_n-1} \int \Big(\frac{p_{j,1}+p_{j,-1}}{2p_{j,0}}\Big)^2p_{j,0}.
\]
 Under $P_{\beta_{\eps}(t)}$  the increments $\Delta_i^n(X,Y-\beta X)$ on block $j$ with $\eps_j=\pm1$ are independent and centered Gaussian with covariance matrix
\[\Sigma_{\eps_j}^i=\begin{pmatrix} \int_{(i-1)/n}^{i/n} \sigma_t^2dt & \eps_j\int_{(i-1)/n}^{i/n} \sigma_t^2\abs{\beta_\eps(t)-\beta}\,dt\\
\eps_j\int_{(i-1)/n}^{i/n} \sigma_t^2\abs{\beta_\eps(t)-\beta}\,dt & \int_{(i-1)/n}^{i/n} (\widetilde{\sigma}_t^2+\sigma_t^2(\beta_\eps(t)-\beta)^2)\,dt
\end{pmatrix},
\]
which implies correlation of $\eps_j\rho_i$. Denoting by $\Sigma_0^i$ the covariance matrix under $H_0$ we obtain therefore
\[ (\Sigma_{\eps_j}^i)^{-1}-(\Sigma_0^i)^{-1}=\frac{-\eps_j\rho_i}{\det(\Sigma_0^i)^{1/2}} \begin{pmatrix} 0&1\\ 1&0\end{pmatrix}.
\]
Using the above, upon denoting with $(Z_{x,i},Z_{y,i})_i$ a sequence of independent standard Gaussian random vectors, we obtain
\begin{align*}
&\int \Big(\frac{p_{j,1}+p_{j,-1}}{2p_{j,0}}\Big)^2p_{j,0}\\
&= E\Big[\Big(\prod_{i=jk_n+1}^{(j+1)k_n} (1-\rho_i^2)\Big)\Big(\frac12\exp\Big(\sum_{i=jk_n+1}^{(j+1)k_n}\rho_iZ_{x,i}Z_{y_i}\Big)
+\frac12\exp\Big(-\sum_{i=jk_n+1}^{(j+1)k_n}\rho_iZ_{x,i}Z_{y_i}\Big) \Big)^2 \Big]\\
&=\Big(\prod_{i=jk_n+1}^{(j+1)k_n} (1-\rho_i^2)\Big)\frac14E\Big[\exp\Big(\sum_{i=jk_n+1}^{(j+1)k_n}2\rho_iZ_{x,i}Z_{y_i}\Big)
+\exp\Big(-\sum_{i=jk_n+1}^{(j+1)k_n}2\rho_iZ_{x,i}Z_{y_i}\Big)+2\Big]\\
&=\Big(\prod_{i=jk_n+1}^{(j+1)k_n} (1-\rho_i^2)\Big)\frac12E\Big[\exp\Big(\sum_{i=jk_n+1}^{(j+1)k_n}2\rho_i^2Z_{y_i}^2\Big)
+1\Big]\\
&=\Big(\prod_{i=jk_n+1}^{(j+1)k_n} (1-\rho_i^2)\Big)\frac12\Big(1+\prod_{i=jk_n+1}^{(j+1)k_n}(1-4\rho_i^2)^{-1/2}\Big)\\
&=\Big(1-\sum_{i} \rho_i^2+\sum_{i\not= j}\rho_i^2\rho_j^2\Big)\Big(1+\sum_i\rho_i^2+3\sum_i\rho_i^4+2\sum_{i\not=j}\rho_i^2\rho_j^2\Big)+O(k_n\max_i\rho_i^6)\\
&=1+2\Big(\sum_{i=jk_n+1}^{(j+1)k_n}\rho_i^2\Big)^2+O(k_n\max_i\rho_i^6),
\end{align*}
where we applied a Taylor expansion to the logarithm of the product, using that $k_n\max_i\rho_i^2=O(nJ_n^{-1-2\alpha})$ is small. Noting $k_n/n\to 0$ and the continuity of the integrands, we have the Riemann sum approximation
\[ \sum_{j=0}^{J_n-1}\Big(\sum_{i=jk_n+1}^{(j+1)k_n}\rho_i^2\Big)^2\approx nk_n\int_0^1
\frac{\sigma_t^4(\beta_\eps(t)-\beta)^4}{(\widetilde{\sigma}_t^2+\sigma_t^2(\beta_\eps(t)-\beta)^2)^2}\,dt.
\]
A similar expansion of the product as above  thus yields the total asymptotic bound
\[ (nk_n)^{-1}\Big(\prod_{j=0}^{J_n-1} \int \Big(\frac{p_{j,1}+p_{j,-1}}{2p_{j,0}}\Big)^2p_{j,0}-1\Big)\to  2\int_0^1\frac{\sigma_t^4(\beta_\eps(t)-\beta)^4} {(\widetilde{\sigma}_t^2+\sigma_t^2(\beta_\eps(t)-\beta)^2)^2}\, dt
\]
for $n,J_n,k_n\to\infty$. Noting $\abs{\beta_\eps(t)-\beta}\le J_n^{-\alpha}\norm{K}_\infty$, the last expression, when scaled up by $nk_n$, is less than $(1-\gamma)^2$ for
$J_n^{-4\alpha}\le C(1-\gamma)^2J_nn^{-2}$ with some constant $C=C(K,\sigma,\widetilde{\sigma})>0$. Hence, for $J_n$ at most $(C(1-\gamma^2))^{-1/(4\alpha+1)}n^{2/(4\alpha+1)}$ the minimax error is bounded by $\gamma$. Choosing $J_n\in\N$ of that order, the separation bound $cJ_n^{-\alpha}\norm{K}_{L^2}$ of the alternative is $\tilde\Gamma r_n$ with $\tilde\Gamma=c(C(1-\gamma)^2)^{-1/(4\alpha+1)}$, as asserted.\qedd

\begin{small}
\bibliographystyle{chicago}
\bibliography{tcb}
\end{small}

\end{document}